\documentclass[]{article}
\usepackage[english]{babel}
\usepackage[T1]{fontenc}
\usepackage[utf8]{inputenc}
\usepackage{graphicx}

\usepackage{tikz}
\usepackage{pgfplots}
\usepackage{pgfplotstable,booktabs,colortbl}
\usetikzlibrary{patterns,positioning,backgrounds,fit}

\usepgfplotslibrary{external}
\tikzset{external/optimize=false}
\tikzset{external/system call={pdflatex
     -shell-escape -halt-on-error
     -interaction=batchmode -jobname "\image" "\texsource"}}
\pgfplotsset{compat=1.14}

\usepackage{expl3}
\ExplSyntaxOn%
\cs_new_eq:NN \fpeval \fp_eval:n%
\ExplSyntaxOff%

\makeatletter
\DeclareRobustCommand*{\escapeus}[1]{%
  \begingroup\@activeescapeus\scantokens{#1\endinput}\endgroup}
\begingroup\lccode`\~=`\_\relax
   \lowercase{\endgroup\def\@activeescapeus{\catcode`\_=\active \let~\_}}
\makeatother
\makeatletter
\DeclareRobustCommand*{\dotus}[1]{%
  \begingroup\@activedotus\scantokens{#1\endinput}\endgroup}
\begingroup\lccode`\~=`\_\relax
   \lowercase{\endgroup\def\@activedotus{\catcode`\_=\active \let~.}}
\makeatother
\makeatletter
\DeclareRobustCommand*{\spaceus}[1]{%
  \begingroup\@activespaceus\scantokens{#1\endinput}\endgroup}
\begingroup\lccode`\~=`\_\relax
   \lowercase{\endgroup\def\@activespaceus{\catcode`\_=\active \let~\ }}
\makeatother

\usepackage[a4paper,top=3cm,bottom=2cm,left=3cm,right=3cm,marginparwidth=1.75cm]{geometry}

\newcommand{\enquote}[1]{``#1''}

\newif\ifloadtikztodos\loadtikztodostrue
\usepackage[small,bf]{caption}
\usepackage{enumitem}

\usepackage{xcolor}

\usepackage{url}

\usepackage[normal,tight]{subfigure}
\makeatletter
  \def\zaptype#1{%
  \listsubcaptions % Finish the last set of subfloats before
  \def\@captype{#1}}% switching to another float type.
\makeatother
\usepackage{amsmath}
\usepackage{amssymb}
\usepackage{amsfonts}

\usepackage{ifthen}
\usepackage[numbers]{natbib}
\usepackage{nameref}

\ifloadtikztodos
  \usepackage[textwidth=2.5cm, textsize=tiny]{todonotes}
    \usepackage{letltxmacro}
    \LetLtxMacro{\oldtodo}{\todo}
    \renewcommand{\todo}[2][]{\tikzexternaldisable\oldtodo[#1]{#2}\tikzexternalenable}
\else
  \usepackage{todo}
  
\fi%
\usepackage{listings}
\makeatletter
\@ifundefined{float@listhead}{}{%
    \renewcommand*{\lstlistoflistings}{%
        \begingroup
    	    \if@twocolumn
                \@restonecoltrue\onecolumn
            \else
                \@restonecolfalse
            \fi
            \float@listhead{\lstlistlistingname}%
            \setlength{\parskip}{\z@}%
            \setlength{\parindent}{\z@}%
            \setlength{\parfillskip}{\z@ \@plus 1fil}%
            \@starttoc{lol}%
            \if@restonecol\twocolumn\fi
        \endgroup
    }%
}
\makeatother
\makeatletter
\def\s@btitle{\relax}
\def\subtitle#1{\gdef\s@btitle{#1}}
\def\@maketitle{%
  \newpage
  \null
  \vskip 2em%
  \begin{center}%
  \let \footnote \thanks
    {\LARGE \@title \par}%
                \if\s@btitle\relax
                \else\typeout{[subtitle]}%
                        \vskip .5pc
                        \begin{large}%
                                \textsl{\s@btitle}%
                                \par
                        \end{large}%
                \fi
    \vskip 1.5em%
    {\large
      \lineskip .5em%
      \begin{tabular}[t]{c}%
        \@author
      \end{tabular}\par}%
    \vskip 1em%
    {\large \@date}%
  \end{center}%
  \par
  \vskip 1.5em}
\makeatother 
\definecolor{DarkGray}{rgb}{0.1,0.1,0.1}
\definecolor{BlueViolet}{rgb}{0.2,0.,0.4}
\definecolor{Grey}{rgb}{0.3,0.3,0.3}
\definecolor{Red}{rgb}{0.8,0.,0.}
\definecolor{Yellow}{rgb}{0.,0.4,0.4}
\definecolor{RoyalBlue}{rgb}{0.,0.14,0.8}
\definecolor{BrickRed}{rgb}{0.8,0.34,0.34}
\definecolor{LimeGreen}{rgb}{0.,0.65,0.25}
\definecolor{Sepia}{rgb}{0.4,0.3,0.3}

\usepackage[acronym,translate=false,toc,counter=equation]{glossaries}
\usepackage{ifnextok}

\def\sym{\sympage}

\glsdisablehyper

\def\glosentryexists#1{%
  \ifglsentryexists{def:#1}{%
    \def\myglossaryprefix{def}%
    \ifglsentryexists{glos:#1}{%
      \glsadd[format=hyperbf,counter=page]{glos:#1}%
    }{%
      \relax%
    }%
  }{%
    \def\myglossaryprefix{glos}%
  }%
}

\makeatletter
\def\glosname#1{%
\glosentryexists{#1}%
\IfNextToken\bgroup%
{\glosnametwoargs@aux{#1}}{\glosnameonearg@only{#1}}%
}%
\def\glosnameonearg@only#1{\glsname[format=hyperbf,counter=page]{\myglossaryprefix:#1}}%
\def\glosnametwoargs@aux#1#2{\glsname[format=hyperbf,counter=page]{\myglossaryprefix:#1}[#2]}%
\makeatother

\makeatletter
\def\glostext#1{%
\glosentryexists{#1}%
\IfNextToken\bgroup%
{\glostexttwoargs@aux{#1}}{\glostextonearg@only{#1}}%
}%
\def\glostextonearg@only#1{\gls[format=hyperbf,counter=page]{\myglossaryprefix:#1}}%
\def\glostexttwoargs@aux#1#2{\gls[format=hyperbf,counter=page]{\myglossaryprefix:#1}[#2]}%
\makeatother

\makeatletter
\def\glospl#1{%
\glosentryexists{#1}%
\IfNextToken\bgroup%
{\glospltwoargs@aux{#1}}{\glosplonearg@only{#1}}%
}%
\def\glosplonearg@only#1{\glspl[format=hyperbf,counter=page]{\myglossaryprefix:#1}}%
\def\glospltwoargs@aux#1#2{\glspl[format=hyperbf,counter=page]{\myglossaryprefix:#1}[#2]}%
\makeatother

\makeatletter
\def\glosfirst#1{%
\glosentryexists{#1}%
\IfNextToken\bgroup%
{\glosfirsttwoargs@aux{#1}}{\glosfirstonearg@only{#1}}%
}%
\def\glosfirstonearg@only#1{\glsfirst[format=hyperbf,counter=page]{\myglossaryprefix:#1}}%
\def\glosfirsttwoargs@aux#1#2{\glsfirst[format=hyperbf,counter=page]{\myglossaryprefix:#1}[#2]}%
\makeatother

\makeatletter
\def\glosfirstpl#1{%
\glosentryexists{#1}%
\IfNextToken\bgroup%
{\glosfirstpltwoargs@aux{#1}}{\glosfirstplonearg@only{#1}}%
}%
\def\glosfirstplonearg@only#1{\glsfirstplural[format=hyperbf,counter=page]{\myglossaryprefix:#1}}%
\def\glosfirstpltwoargs@aux#1#2{\glsfirstplural[format=hyperbf,counter=page]{\myglossaryprefix:#1}[#2]}%
\makeatother

\makeatletter
\def\Glosname#1{%
\glosentryexists{#1}%
\IfNextToken\bgroup%
{\Glosnametwoargs@aux{#1}}{\Glosnameonearg@only{#1}}%
}%
\def\Glosnameonearg@only#1{\Gls[format=hyperbf,counter=page]{\myglossaryprefix:#1}}%
\def\Glosnametwoargs@aux#1#2{\Gls[format=hyperbf,counter=page]{\myglossaryprefix:#1}[#2]}%
\makeatother

\makeatletter
\def\Glospl#1{%
\glosentryexists{#1}%
\IfNextToken\bgroup%
{\Glospltwoargs@aux{#1}}{\Glosplonearg@only{#1}}%
}%
\def\Glosplonearg@only#1{\Glspl[format=hyperbf,counter=page]{\myglossaryprefix:#1}}%
\def\Glospltwoargs@aux#1#2{\Glspl[format=hyperbf,counter=page]{\myglossaryprefix:#1}[#2]}
\makeatother

\makeatletter
\def\GLOSname#1{%
\glosentryexists{#1}%
\IfNextToken\bgroup%
{\GLOSnametwoargs@aux{#1}}{\GLOSnameonearg@only{#1}}%
}%
\def\GLOSnameonearg@only#1{\GLS[format=hyperbf,counter=page]{\myglossaryprefix:#1}}%
\def\GLOSnametwoargs@aux#1#2{\GLS[format=hyperbf,counter=page]{\myglossaryprefix:#1}[#2]}%
\makeatother

\makeatletter
\def\GLOSpl#1{%
\glosentryexists{#1}%
\IfNextToken\bgroup%
{\GLOSpltwoargs@aux{#1}}{\GLOSplonearg@only{#1}}%
}%
\def\GLOSplonearg@only#1{\GLSpl[format=hyperbf,counter=page]{\myglossaryprefix:#1}}%
\def\GLOSpltwoargs@aux#1#2{\GLSpl[format=hyperbf,counter=page]{\myglossaryprefix:#1}[#2]}%
\makeatother


\def\dummysympage{\sympage}
\makeatletter
\def\glosformula#1{%
  \glosentryexists{#1}%
  \ifx\sym\dummysympage%
    \glosformula@nonnumbered{#1}%
  \else%
    \glosformula@numbered{#1}%
  \fi%
}%
\def\glosformula@numbered#1{%
\glsadd[format=hyperit,counter=equation]{\myglossaryprefix:#1}%
\glsentryuseri{\myglossaryprefix:#1}%
}%
\def\glosformula@nonnumbered#1{%
\glsadd[format=hyperbf,counter=page]{\myglossaryprefix:#1}%
\glsentryuseri{\myglossaryprefix:#1}%
}%
\makeatother

\newenvironment{myequation*}{\def\sym{\sympage}\equation\aligned\protect}{\nonumber\endaligned\endequation}%
\makeatletter
\newcommand\mynonumber{\start@align\@ne\st@rredtrue\m@ne}%
\newcommand\mynumber{\start@align\@ne\st@rredfalse\m@ne}%
\makeatother

\newenvironment{myalign*}{\def\sym{\sympage}\mynonumber\protect}{\endalign}%
\newcommand{\symeq}[1]{
  \protect\gls[format=hyperit]{sym:#1}
  \ifglsentryexists{def:#1}{%
    \glsadd[format=hyperbf,counter=page]{def:#1}
  }{%
    \relax%
  }
}
\newcommand{\sympage}[1]{\glssymbol[format=hyperbf,counter=page]{def:#1}}

\usepackage{dsfont}
\def\unity{\mathds{1}}
\def\reel{\mathbb{R}}

\def\citesection{sect.}
\def\citepage{p.}

\DeclareMathOperator{\var}{var}
\DeclareMathOperator{\cov}{cov}

\def\vp{\varphi}
\def\n{\nabla}
\def\de{\Delta}
\newcommand{\R}{\mathbb{R}}
\newcommand{\E}{\mathbb{E}}

\newcommand{\Lop}{\mathcal{L}}

\def\TF{\emph{TensorFlow}} % TensorFlow framework
\def\Pyth{\emph{Python}}% Python language
\newcommand{\python}[1]{\texttt{#1}}% for python objects, classes, ...
\newcommand{\module}[1]{\texttt{#1}}% for python modules
\newcommand{\num}[1]{\pgfmathprintnumber[fixed,int detect]{#1}{}}% for printing numbers with good readibility

\newglossaryentry{glos:BBGD}{
	type=\acronymtype,
	name={Barzilei-Borwein Gradient Descent},
	text={BBGD},
	first={Barzilei-Borwein Gradient Descent (BBGD)},
	description={},
	sort=BBGD
}

\newglossaryentry{glos:CLT}{
	type=\acronymtype,
	name={Central Limit Theorem},
	text={CLT},
	first={Central Limit Theorem (CLT)},
	description={},
	sort=CLT
}

\newglossaryentry{glos:EQN}{
	type=\acronymtype,
	name={ensemble quasi-Newton},
	text={EQN},
	first={ensemble quasi-Newton (EQN)},
	description={},
	sort=EQN
}

\newglossaryentry{glos:GD}{
	type=\acronymtype,
	name={Gradient Descent},
	text={GD},
	first={Gradient Descent (GD)},
	description={},
	sort=GD
}

\newglossaryentry{glos:GLA}{
	type=\acronymtype,
	name={Geometric Langevin Algorithm},
	text={GLA},
	first={Geometric Langevin Algorithm (GLA)},
	description={},
	sort=GLA
}

\newglossaryentry{glos:HMC}{
	type=\acronymtype,
	name={Hamiltonian Monte Carlo},
	text={HMC},
	first={Hamiltonian Monte Carlo (HMC)},
	description={},
	sort=HMC
}

\newglossaryentry{glos:IAT}{
	type=\acronymtype,
	name={Integrated Autocorrelation Time},
	text={IAT},
	first={Integrated Autocorrelation Time (IAT)},
	description={},
	sort=IAT
}

\newglossaryentry{glos:ISST}{
	type=\acronymtype,
	name={Infinite Switch Simulated Tempering},
	text={ISST},
	first={Infinite Switch Simulated Tempering (ISST)},
	description={},
	sort=ISST
}

\newglossaryentry{glos:MCMC}{
	type=\acronymtype,
	name={Markov Chain Monte Carlo},
	text={MCMC},
	first={Markov Chain Monte Carlo (MCMC)},
	description={},
	sort=MCMC
}

\newglossaryentry{glos:MH}{
	type=\acronymtype,
	name={Metropolis-Hastings},
	text={MH},
	first={Metropolis-Hastings (MH)},
	description={},
	sort=MH
}

\newglossaryentry{glos:PES}{
	type=\acronymtype,
	name={Potential Energy Surface},
	text={PES},
	first={Potential Energy Surface (PES)},
	description={},
	sort=PES
}

\newglossaryentry{glos:QMC}{
	type=\acronymtype,
	name={Quasi Monte Carlo},
	text={QMC},
	first={Quasi Monte Carlo (QMC)},
	description={},
	sort=QMC
}

\newglossaryentry{glos:SD}{
	type=\acronymtype,
	name={Steepest Descent},
	text={SD},
	first={Steepest Descent (SD)},
	description={},
	sort=SD
}

\newglossaryentry{glos:SGD}{
	type=\acronymtype,
	name={Stochastic Gradient Descent},
	text={SGD},
	first={Stochastic Gradient Descent (SGD)},
	description={},
	sort=SGD
}

\newglossaryentry{glos:SGHMC}{
	type=\acronymtype,
	name={Stochastic Gradient Hamiltonian Monte Carlo},
	text={SGHMC},
	first={Stochastic Gradient Hamiltonian Monte Carlo (SGHMC)},
	description={},
	sort=SGHMC
}

\newglossaryentry{glos:SGLD}{
	type=\acronymtype,
	name={Stochastic Gradient Langevin Dynamics},
	text={SGLD},
	first={Stochastic Gradient Langevin Dynamics (SGLD)},
	description={},
	sort=SGLD
}

\newglossaryentry{glos:TATi}{
	type=\acronymtype,
	name={Thermodynamic Analytics ToolkIt},
	text={TATi},
	first={Thermodynamic Analytics ToolkIt (TATi)},
	description={},
	sort=TATi
}

\makeglossaries

\newcommand{\myexperimentspath}{/nosave/heber/Experiments/ThermodynamicAnalyticsToolkit}%

\title{TATi-Thermodynamic Analytics ToolkIt: TensorFlow-based software for posterior sampling in machine learning applications}
\author{Frederik Heber, Zofia Trstanova, Benedict Leimkuhler}

\begin{document}

\maketitle

\begin{abstract}
With the advent of GPU-assisted hardware and maturing high-efficiency software platforms such as TensorFlow and PyTorch, Bayesian posterior sampling for neural networks becomes plausible.
In this article  we discuss Bayesian parametrization in machine learning based on Markov Chain Monte Carlo methods, specifically discretized stochastic differential equations such as Langevin dynamics and extended system methods in which an ensemble of walkers is employed to enhance sampling. We provide a glimpse of the potential of the sampling-intensive approach by  studying (and visualizing) the loss landscape of a neural network applied to the MNIST data set. Moreover, we investigate how the sampling efficiency itself can be significantly enhanced through an ensemble quasi-Newton preconditioning method.  This article accompanies the release of a new TensorFlow software package, the Thermodynamic Analytics ToolkIt, which is used in the computational experiments.
\end{abstract}

\section{Introduction}
% The problem setting:
% Neural networks are employed everywhere, trained through optimization,
% while it is not understood why the non-convex optimization actually works.
%
%With this article we introduce a new software library, the Thermodynamic Analytics ToolkIt (TATi) which interfaces to \TF{} to allow study of posterior distributions of statistical models\todo{fh: This statement is too general and thereby applies exactly to what TFP offers.}.  As primary application, we focus on supervised learning in neural networks simply because such problems are the most directly accessible and easiest to describe of \TF{}'s targets, but most of what is said throughout this article could be extended to a significantly larger class of parameterized models used in machine learning applications.

The fundamental role of neural networks (NNs) is readily apparent from their
widespread use in machine learning in applications such as natural language processing~(\cite{natlang}), social network analysis~(\cite{socnets}), medical
diagnosis~(\cite{autdiag,heartdiag}), vision systems~(\cite{vision}), and robotic path planning~(\cite{lecun2015deep}). The greatest success of these models lies in their flexibility, their ability to represent complex, nonlinear relationships in high-dimensional data sets, and the availability of frameworks that allow NNs to  be implemented on rapidly evolving GPU platforms~(\cite{krizhevsky2012imagenet,he2016deep}). The industrial appetite for deep learning has led to very rapid expansion of the subject in recent years, although, as pointed out by ~\cite{dunson2018statistics},  at times the mathematical and theoretical understanding  of these methods has been swept aside in the rush to advance the methodology.   

The potential impact on society of machine learning algorithms demands that their exposition and use be subject to the highest standards of clarity, ease of interpretation, and uncertainty quantification.  Typical NN training seeks 
%specific architecture of neural networks allows to
to optimize the parameters of the network (biases and weights) under the constraint that the training data set is well approximated~(\cite{hardt2015train,Goodfellow2014a}).       In the Bayesian setting, the parameters of a neural network are defined by the observations, but only in the probabilistic sense, thus specific parameter values are only realized as modes or means of the associated distribution, which can require substantial computation. 
Bayesian approaches based on exploration of the posterior probability distribution have been discussed throughout the development of neural networks (\cite{MacKay1992,Neal1992,Hinton1993,Barber1998}), and underpin much of the work in this field,  but they are less commonly implemented in practice for  ``big data'' applications due to legitimate concerns about efficiency (\cite{Welling2011}). While the idea of sampling (or partially sampling) the posterior of large scale neural networks is not new, improvements in computers continually render this goal more plausible.  For a recent discussion, see the PhD thesis of \cite{Gal2016}, which again champions the use of a (Bayesian) statistical framework, mentioning among other aims the prospect for meaningful uncertainty quantification in deep neural networks.   Posterior sampling has the potential for broad impact in several research areas related to NN construction, including: (i) the relationship between network architecture and parameterization efficiency (\citet{Safran2016,Livni2014,Haeffele2017,Soltanolkotabi2017}), (ii) the
visualization the loss manifold and the parameterization process (see \citet{Draxler2018,Im2016,Li2017,Goodfellow2014a}), (iii)
and the assessment of the generalization capability of networks (see \citet{Dinh2017,Hochreiter1997a,Kawaguchi2017,Hoffer2017}).

We have recently developed  the \glostext{TATi}, a python framework whose purpose is to facilitate the sampling of the posterior parameter distribution of automated machine learning systems with a balance of ease of use and computational efficiency, by leveraging highly optimised computational procedures within \TF{}.\footnote{Installation of TATI is as simple as {\tt pip install tati} from the unix command line; the software includes a readable user guide and programmer's manual.}  Although still at an early stage in its development, this software package makes possible the exploration of sampling-based approaches in the high-dimensional setting.  This article provides the motivation for TATi by illustrating that it facilitates analysis of the loss manifold in a way which would be impossible using standard optimization methods.

%In this setting, the accuracy of the model is assessed by estimation of the {\em generalization error}, which is the difference between the {\em empirical loss} (the loss computed over a finite dataset) and the {\em population loss} (the expected loss with respect to the entire data distribution). In other words, the generalization error estimates how accurate the trained model is on  unseen data, given the fitted parameters.    The concept of generalization error illustrates the fundamental statistical nature of the parameterization problem.   
%In order to provide probabilistic guarantees on the predictions of neural networks we adopt a Bayesian perspective and sample from the posterior distribution of the parameters given a training dataset~. 

The standard approach to Bayesian sampling in high dimensions relies on  Markov Chain Monte Carlo methods, but these can be difficult to scale to  large system size.   
%Performing Bayesian inference over the neural network parameters has the potential for better generalization 
%by reducing {\em overfitting} to the training data set~\cite{Welling2011}. {\color{blue} we need a bit more references for posterior sampling} 
 In the context of deep learning, we use the term {\em large}  to refer  both to large data set size (which translates into expensive likelihood computations) and to large numbers of parameters (usually the consequence of adopting a {\em deep learning} paradigm).  We base our software on discretized stochastic differential equations which offer a reliable and accurate means of computing \glostext{MCMC} sampling paths while providing rigorous results with statistical error bounds (\cite{leimkuhler2015computation}).  An underlying assumption in using software like this is that the exploration is limited to a bounded region of the parameter space by some structural features of the problem and/or its regularization.  The methodology we describe could also be extended  to include an explicit localization scheme to implement a quasi-stationary (spatially restricted) distribution \cite{qsd1}, although we do not discuss this here.
%The posterior distribution of a parameter vector $\theta\in \R^{N}$ given a dataset $X$ is
%\[
%\pi(\theta \mid X) \propto \pi(X \mid \theta) \pi(\theta), 
%\]
%where $\pi(X \mid \theta)$ is %the likelihood and $ %\pi(\theta)$ is the prior %distribution. 

As an indication of the possible scope for practical posterior sampling in large scale  machine learning, we note that Markov Chain Monte Carlo methods like those implemented in \glostext{TATi} have been successfully deployed in the more mature setting of very large scale molecular simulations in computational sciences for applications in physics, chemistry, material science and biology, with variables numbering in the millions or even billions (\cite{Klein,namd,Kadau}).  The goal in molecular dynamics is similarly the sampling of a target probability distribution, although the distribution typically arises from semi-empirical modelling of interactions among atoms of the substances of interest. The scientific communities in molecular sciences are developing highly efficient ``enhanced sampling'' procedures to tackle the computational difficulties of large scale sampling (see the surveys (\cite{mds1,mds2}) for some examples of the wide variety of methods being used in biomolecular applications).
In analogy with molecular dynamics, the focus on posterior sampling makes possible the elucidation of reduced descriptions of neural networks through concepts such as free energy calculation (\cite{lelievre}) and transition path sampling (\cite{tps}).  Although we reserve detailed study of these ideas for future work, the fundamental tool in their construction and practical implementation is the ability to compute sample sequences efficiently and reliably using \glostext{MCMC} paths.   Further potential benefits of the posterior sampling approach lie in the fact that it offers a means of high-dimensional uncertainty quantification through standard statistical methodology (\cite{UQ,Gal2016}).

%{\color{blue}In developing TATi, one goal we have had in mind is to improve the graphical analysis of the loss landscape so as to better understand its structure vis a vis the performance of parameterization algorithms.  } 

We favor schemes based on underdamped Langevin dynamics (with canonical momentum variables). Such methods have excellent properties in terms of accuracy of statistical averages and convergence rates~(\cite{leimkuhler2015computation,lelievre}). A variety of other methods are available which are also based on second order dynamics ~(\cite{
chen2014stochastic,patterson2013stochastic,shang2015covariance,matthews2018langevin}).  One of the purposes of the TATi software is to provide a relatively simple mechanism for the implementation and evaluation  of sampling strategies of statistical physics in machine learning applications. \footnote{
While \TF{}'s graph programming structure is well explained and motivated by its developers and provides efficient execution on GPUs, it is not straightforward for numerical analysts, statisticians and statistical physicists to modify in order to test their methods.    \glostext{TATi} therefore uses a simplified interface structure that allows algorithms to be coded directly in pure \Pyth{} and then linked to \TF{} for efficient calculation of gradients for arbitrary \TF{} models.   As a consequence, building and testing a new sampling scheme in \glostext{TATi} requires little knowledge of the underpinnings of \TF{}.}

Another goal we have had in mind in this study is to gain insight into  the loss landscape (given by the log of the posterior density) so as to better understand its structure vis a vis the performance of parameterization algorithms.
The loss is not convex in general~(\cite[\citesection~4]{LeCun1998}), and its corrugated (``metastable'') structure has been compared to models of  spin-glasses ~(\cite{Choromanska2015}):  there is a band of minima close to the global minimum as lower bound and whose multitude diminishes exponentially for larger loss values.  The availability of an efficient sampling scheme gives a means of better understanding the loss landscape and its relation to the network (and properties of the data set). We illustrate some of the potential for such studies in Section 3 of this article, where we examine the loss landscape of an MNIST classification problem.

Finally, we note that the statistical perspective underpins many optimization schemes in current use in machine learning, such as the \glostext{SGD} method. The stochasticity enters into this method through the subsampling of the dataset at every iteration of the optimization algorithm, due to  `minibatching'. Several stochastic optimisation methods have been proposed which aim to improve the computational efficiency or reduce 
generalisation error, e.\,g., RMSprop~(\cite{Hinton2014}), AdaGrad~(\cite{Duchi2011}), Adam~(\cite{kingma2014adam}), entropy-SGD~(\cite{chaudhari2016entropy}). \glostext{SGLD}~(\cite{Welling2011}), and these, similarly, can be given a statistical (sampling) interpretation.  Indeed, when its stepsize is held fixed and under simplifying assumptions on the character of the gradient noise,  it can be viewed as a first-order discretization of overdamped Langevin dynamics.  We discuss \glostext{SGD} and \glostext{SGLD} in Sec 2, in order to motivate SDE sampling methods.

% Enhanced sampling and optimisation schemes: Entropy-SGD~\cite{chaudhari2016entropy} is an optimization method for deep learning designed to exploit the local geometric properties of the objective. It uses local sampling to optimize the stochastig gradient descent with respect to 'local free energy' obtained by averaging the gradients over several Langevin steps.

%\subsection{Towards large scale inference in TensorFlow.}
%The main applications of this approach as 1) bringing more insight into free energy loss lanscapes of neural networks, 2) finding parametrisations with minimal generalisation error, 3) allowing for sensibility analysis for neural networks. We emphasize the use of Langevin dynamics samplers, and enhanced sampling schemes. 
%TATi is an open-source library and aims at connecting the machine learning and the statistical physics community.  

The remainder of this paper is structured as follows: Section~\ref{sec:sampling} provides an overview of \glostext{MCMC},  Langevin dynamics schemes and other sampling strategies and the concepts from numerical error analysis that underpin our approach.  We also outline there the Ensemble Quasi-Newton method (\cite{Matthews2018}). 
Subsequently, in Section 3, we discuss the accuracy of various schemes as well as their convergence behavior in the context of neural network posterior sampling.   We also elaborate on how the approach we describe can be used to handle moderate-dimensional sampling on the MNIST dataset.  The last part of this section consists of a demonstration of the convergence acceleration of the EQN sampler in the MNIST application.    The numerical results presented in the paper represent a proof-of-concept for the utility of the posterior sampling paradigm and the TATi software which implements it.

\section{Markov Chain Monte Carlo methods and stochastic differential equations}\label{sec:sampling}

In this section, we provide a brief introduction to the \glostext{MCMC} methods we have used for sampling high dimensional distributions.  These methods include schemes based on discretization of stochastic differential equations, especially Langevin dynamics. We discuss in some detail the construction of numerical schemes in order to control the finite time stepsize bias.  Moreover, we also describe an ensemble quasi-Newton method implemented in TATi that adaptively rescales the dynamics to enhance sampling of poorly conditioned target posteriors.
%\todo{Note the following notations which need to be consistent in subsequent sections: stepsize is $\varepsilon$, dimension of parameter vector is $N$, parameters live in $\Omega\subset R^N$, dimension of data set is $M$, subsample size $m$, the timestep index is now $k$ instead of $n$.  Also I use subscript $k$ to indicate the indexing of timesteps. $\Omega$ is the parameter space and $\Omega_{\rm ext}$ is the extended parameter space with the addition of momenta.}

Assume that we are given a dataset ${\cal D} = ({\cal X}, {\cal Y}) = \{(x_i,y_i)\}_{i=1}^M$ comprising inputs and outputs of an unknown functional relationship.  We also assume we are given a neural network defined by a vector of parameters $\theta \in \Omega \subset \reel^N$ which acts on an input $x\in {\cal X}$ to produce output $f(x,\theta)$.    The goal of training the network is then typically formulated as solving an optimization problem over the parameters given the dataset:
\begin{equation}\label{math:neuralnetwork-minimization}
	\min_{\theta \in \reel^N} L(\theta , \mathcal{D}),
\end{equation}
where the function $L(\theta, \mathcal{D})$ is the total loss function associated to the dataset ${\cal D}$, defined by
\begin{equation}\label{math:lossfunction}
L(\theta, \mathcal{D})=  \sum_{(x,y) \in {\cal D}} l(f(x, \theta),y) = \sum_{i=1}^M l(f(x_i, \theta),y_i).
\end{equation}
The loss function $l(\hat{y},y)$ depends on the metric choice, for example the squared error, logarithmic loss,  cross entropy loss, etc~(\cite{lossfunctions}).  The total loss is in general not a convex function of the parameters $\theta$ even though $l$ may be convex as a function of $y,\hat{y}$. In general the loss landscape is rough there are likely to be many local minima and flattened intermediate states in the loss manifold $L(\theta, \mathcal{D})$, as well as many saddles, see \cite{Baldi1989,Choromanska2015}. % for an investigation of their structure in the context of multi-layer perceptrons using absolute loss and hinge loss.

The parametrization procedure is based on optimization algorithms, which generate a sequence of  parameter vectors $\theta_0,\theta_1,\ldots, \theta_k,\ldots$, converging to a (local) minimizer of~\eqref{math:neuralnetwork-minimization} as $k \rightarrow \infty$. The basic optimization algorithm is the \glostext{GD}, which uses the negative gradient to update the parameter values, i.\,e.,  
\begin{equation}
\theta_{k+1} =  \theta_k - \varepsilon \nabla L(\theta_k),
\label{eq:gd_update_step}
\end{equation}
 where $\varepsilon$ is the stepsize (or learning rate), which is either constant or may be varied during the computation. \glostext{GD} converges for a convex function and for a smooth non-convex total loss function it converges to the nearest local minimum~(\cite{Nocedal1999}).

Gradient calculations are the primary computational burden when training neural networks, i.\,e., the computational cost of a gradient can be taken as a reasonable measure of computational work.  As the total loss $L$ implicitly depends on the whole dataset, one natural idea that reduces the computational cost is to exploit the redundancy in the dataset by estimating the gradient of the average loss from a subset of the data, that is, to replace the gradient in each parameterization step by the approximation
\[
\nabla\widetilde{L}(\theta)\approx \tfrac M m \sum_{i\in S_k} \nabla l(f(x_i, \theta), y_i).
\]
Where $S_k$ represents a randomized data subset (re-randomized throughout the training process) of dimension $m$, this method, which has many variants, is referred to as \glostext{SGD}. 

We are interested in the Bayesian inference formulation, where $\theta$ is the parameter vector and we wish to sample from the posterior distribution $\pi(\theta \mid {\cal D})$ of the parameters given a dataset of size $M>0$,
\[
\pi(\theta \mid {\cal D}) \propto \pi_0(\theta)\prod_{i=1}^M \pi( (x_i,y_i) \mid \theta),
\]
with prior probability density $\pi_0(\theta)$, and likelihood $\pi((x,y)\mid \theta)$. 

\glosfirst{SGLD} (\cite{Welling2011}) generates a step based on a subset of the data and injects additional noise. 
The additive noise creates a controllable stochastic model with known ergodic properties and it improves the numerical stability and convergence properties of the training algorithm.  For additional discussion of these points see \cite{LeMaVl2019}, where small injected noise has been shown to substantially accelerate and improve robustness of the training process for neural network models.

The \glostext{SGLD} parameter update is using a sequence of stepsizes $\{\varepsilon_n\}$ and reads 
\begin{equation}\label{eq:sgld}
\theta_{k+1} = \theta_{k} - \varepsilon_k\nabla\tilde{L}(\theta_k)  + \sqrt{2\beta^{-1} \varepsilon_k} G_k, 
\end{equation}
with 
\[
 \nabla\tilde{L}(\theta) = - \beta^{-1}   \nabla \log \bigl (\pi_0(\theta) \bigr ) - \beta^{-1}  \frac{M}{m}\prod_{i\in S_k}\nabla \log \bigl (\pi((x_i,y_i) \mid \theta) \bigr ) 
\]
where $\beta>0$ is a constant (in physics, it would be associated with reciprocal temperature) and again $S_k$ is a random subset of indices ${1, \ldots, M}$ of size $m$. 
   The original algorithm uses a diminishing stepsize sequence, however, \cite{vollmer2016exploration} showed that fixing the stepsize has the same efficiency, up to a constant.   
   
   %Experiments in~\cite{Neelakantan2015} demonstrate the effectiveness of this method in avoiding trapping in local minima.
   
%		As \glsentrytext{glos:SGD} is realized by an ML or MAP estimation it tends to overfit. As a remedy in \cite{Welling2011} it was proposed to add another noise term that dominates the noise injected through the stochastic gradient calculation. This allows to prevent a collapse of the solution to the MAP estimate. This method is called , in \cite{Teh2016} its convergent properties are described. Furthermore, in \cite{Neelakantan2015} there are many numerical experiments that highlight its effectiveness in overcoming being trapped in local minima.
%		The extent by which these stochastic gradient methods are able to escape the nearest local minima is determined by the magnitude of the injected noise. However, this fully depends on the dataset and on the size of the subset. Therefore, it can not be controlled sensibly.
%		It may be argued whether the global minimum of some training dataset is actually desirable as it most likely is drastically overfit. However, undoubtedly some local minimum is just as undesirable. This is especially true in the current situation where accountability of neural networks is highly desired. Eventually, a global optimization is required. Such an approach is realized through sampling of the loss manifold using proper dynamics such as Langevin or Hamiltonian that maintain the canonical distribution. %%%%%%%%%%%%%%%%%%%%%%%%%%%%%%%%%%%%%%%%%%

Under the assumption that the gradient noise is uncorrelated and identically distributed from step to step, it is easy to demonstrate that \glostext{SGLD} with fixed stepsize is an Euler-Maruyama (first order) discretization of overdamped Langevin dynamics (\cite{Welling2011}):
\begin{equation}
d \theta=-  \nabla\tilde{L}(\theta)dt+\sqrt{2\beta^{-1}}dW_t.
\label{eq: overdamped Langevin}
\end{equation}
This creates a natural starting point for our approach: we simply change the perspective to focus on the sampling of the posterior distribution, rather than the identification of its mode.

In case a full gradient is used, dynamics~\eqref{eq: overdamped Langevin} preserves the Boltzmann distribution with measure proportional to an exponential function of the negative loss:
\begin{equation}
\pi(d\theta) \propto {\rm e}^{-\beta L(\theta)}d\theta.
\label{boltzmann}
\end{equation}

Given a generic process providing samples $\theta_k$ asymptotically distributed with respect to a defined target measure $\pi$, it is possible to use sampling paths to estimate integrals with respect to $\pi$. Define the finite time average of a $C^{\infty}$ function $\varphi$ by
\[
\widehat{\varphi}(K) := \frac{1}{K+1}\sum_{k=0}^{K} \varphi\left(\theta_k\right).
\]
For an ergodic process, we have
	\begin{equation}
	\lim_{K\rightarrow \infty} \widehat{\varphi}(K)=\int_{\Omega}\varphi(\theta) \, \pi(d \theta).
	\label{eq: ergodic averages markov chain}
	\end{equation}

The convergence rate of the limit above is given by the \glostext{CLT}. Given a generic process generating samples $\theta_k$ from a target distribution with density $\pi$, the variance of an observable $\varphi$ behaves, asymptotically for large $K$, as 
\[
\var\widehat{\varphi}(K) \sim
{\tau}_{\varphi}\var{\varphi}/K
\]
 where $\tau_{\varphi}$ is the {\em integrated autocorrelation time}, see \citet[\citesection~3]{Goodman2010}.  $\tau_{\varphi}$ can be viewed as a measure of the redundancy of the sampled values or the number of steps until the sampled values of $\varphi$ decorrelate, thus $\tau_{\varphi} = 1$ is optimal, i.\,e., immediately stepping from one independent state to the next. The \glostext{IAT} $\tau_{\varphi}$ can also be calculated as 
\begin{equation}
	\label{math:integrated_autocorrelation_time}
	\tau_{\varphi} = 1+ 2\sum^{\infty}_{1} \frac{ C_{\varphi}(k)} {C_{\varphi}(0)} \quad	\text{with} \quad C_{\varphi}(k) =  \cov[\varphi \bigl (\theta_{k} \bigr ), \varphi \bigl (\theta_{0} \bigr)],
\end{equation}
where the covariance is averaged over the initial condition.

In practice, the process of discretization of the SDE introduces asymptotic bias, which is controlled by the stepsize, however it is nonetheless possible and useful to compute the IAT in such a case to describe the convergence to the asymptotically perturbed equilibrium distribution.

%The above holds also for sampling approaches based on Langevin Dynamics. 
%There, we may use the \glostext{IAT} to gauge the \emph{exploration speed} for each sampled trajectory $X(t)$.

\paragraph{Langevin Dynamics.}	

Sampling from~\eqref{boltzmann} is one of the main challenges of computational statistical physics.  
%The problem is known to be very difficult when the dimensionality is high and the distribution multimodal (a consequence of the loss function $L$ having multiple local minima).   
There are three popular alternative approaches to sample from~\eqref{boltzmann}:  discretization of continuous stochastic differential equations (SDEs), Metropolis-Hastings based algorithms and deterministic dynamics, as well as combinations among the three groups. 
%In the \glostext{TATi} sofware, we have implemented methods derived from SDEs such as \eqref{eq: overdamped Langevin}, underdamped Langevin schemes, but also Hamiltonian Monte Carlo methods.

Langevin dynamics is an extended version of~\eqref{eq: overdamped Langevin}:
\begin{equation}
  \label{eq: Langevin}
  \left\{
  \begin{aligned}
  d\theta_t & =  M^{-1}p_t \, dt, \\
  dp_t & = -\n L(\theta_t) \, dt - \gamma  M^{-1}p_t \, dt + \sigma \, dW_t,
  \end{aligned}
\right.
\end{equation}
where $p$ is the momentum variable and $\gamma>0$ is the friction. The fluctuation-dissipation relation $\sigma^2={\frac{2\gamma}{\beta}}$ ensures that the extended (canonical) distribution with density
\[
\pi_{\rm ext} \propto \exp \left [ -\beta \left ( p^TM^{-1}p/2  \right )\right ]
\pi(\theta)
\]
is preserved; the target distribution is recovered by marginalization.   

Whereas the momentum is a physical variable in statistical mechanics, it is introduced as an artificial auxiliary variable in the machine learning application. $M$ in \eqref{eq: Langevin} is a positive definite  mass matrix, which can in many cases be taken to the identity matrix. 

  %Many alternative definitions of the momenta could be implemented \cite{alternativemomenta}, although in the current version of \glostext{TATi} the elementary definition is used.

The discretization of stochastic dynamics introduces bias in the invariant distribution (see below, for discussion).  Although the bias can be removed through the incorporation of a  \glostext{MH} step, such  methods do not always scale well with the dimension~\cite{Beskos2009,Robert2015}; and the \glostext{MH} test is often neglected in practice.  For example the \enquote{unadjusted Langevin algorithm} of \cite{Moulines} tolerates the presence of small stepsize-dependent bias, in order to obtain faster convergence (reduced asymptotic variance or lower integrated autocorrelation time for observables of interest) and better overall efficiency. 

%In \glostext{TATi}, we emphasize Langevin dynamics discretizations but \glostext{TATi} includes a Hamiltonian Monte Carlo method (see below) which we briefly test in this article; the Metropolis scheme implemented in \glostext{TATi}'s HMC could be adapted to other, more general proposal distributions.

\paragraph{Systematic design of schemes.}

The mathematical foundation for the construction of discretization schemes for~\eqref{eq: Langevin} by splitting of the generator of the dynamics is now well understood~(\cite{leimkuhler2015computation}). 
Although different choices can be made, the usual starting point is an additive decomposition of the generator of dynamics~\eqref{eq: Langevin} into three operators $\Lop = A + B + O$, where
\begin{equation}
A:=M^{-1}p \cdot \n_\theta, \qquad B:=-\n  L(\theta)\cdot \n_p, \qquad O:= -\gamma M^{-1}p\cdot \n_p + \sigma \de_p \,,
\label{eq: A B O operators}
\end{equation}
The main idea is that each of these sub-dynamics can be resolved exactly in the weak (distributional) sense. Note that the (\enquote{O} step) dynamics
\begin{equation}
dp_t  =- \gamma M^{-1}p_t \, dt + \sigma \, dW_t,
\label{eq: FD part}
\end{equation}
has an analytical solution
\begin{equation}
p_t= \alpha_{ t}{p}_0+ \sigma \int_0^t\alpha_{t-s} dW_s\,,\quad \alpha_{t}:={\rm e}^{-\gamma M^{-1} t}.
\label{eq: analytical solution fd}
\end{equation}
A Lie-Trotter splitting of the elementary evolution generated by $A,B,$ and $O$ provides six possible first-order splitting schemes of the general form
\[
P^{Z,Y,X}_{\varepsilon} = {\rm e}^{\varepsilon Z} \,{\rm e}^{\varepsilon Y} \,{\rm e}^{\varepsilon X},
\]
with all possible permutations $(Z,Y,X)$ of $(A, B, O)$, and second-order splitting schemes are then obtained by a Strang splitting of the elementary evolutions generated by $A,B,$ and $O$. %There are six possible schemes in both cases.

Using the notation for the three operators of the sub-dynamics~\eqref{eq: A B O operators}, we define the following updates which will be combined in the full scheme for the discretization of the Langevin dynamics with stepsize $\varepsilon$:
\begin{equation}
\label{eq:ABO_steps}
\begin{aligned}
& A_{\varepsilon}: \theta \rightarrow \theta + \varepsilon p , \\
& B_{\varepsilon}: p \rightarrow p -\varepsilon \nabla L(\theta), \\
& O_{\varepsilon} : p \rightarrow \alpha p + \sqrt{\beta^{-1}(1-\alpha^2)}R,
\end{aligned}
\end{equation}
with $\alpha={\rm e}^{-\gamma \varepsilon}$ and $R\sim\mathcal{N}(0,1)$. A numerical methods is easily specified by a string such as `ABO'.  This is an instance of the Geometric Langevin Algorithm (\cite{bou2010long}).    We also consider symmetric compositions of several basic steps.
The `BAOAB' scheme is in this notation $B_{\varepsilon/2}A_{\varepsilon/2}O_{\varepsilon}A_{\varepsilon/2}B_{\varepsilon/2}$.   
This method can be written out in detail as a step from $(\theta_k, p_k)$ to $(\theta_{k+1},p_{k+1})$ as follows:

\begin{eqnarray*}
p_{k+1/2} & = &p_k - \tfrac{\varepsilon}2 \nabla L(\theta_k),\\
\theta_{k+1/2} & = &\theta_{k} + \tfrac{\varepsilon}2 M^{-1}p_{k+1/2},\\
p_{k+1} & = & \alpha p_{k+1/2}+\sqrt{\beta^{-1} (1-\alpha^2)}R_{k+1/2},\\
\theta_{k+1} & = &\theta_{k+1/2} + \tfrac{\varepsilon}2 M^{-1}p_{k+1},\\
p_{k+1} & = &p_{k+1/2} - \tfrac{\varepsilon}2\nabla L(\theta_{k+1}).
\end{eqnarray*}

Given the sequence of samples $\left(\theta_k, p_k\right)$ determined using such a method, we can  approximate expected values of a $C^{\infty}$ function of $\theta$ and $p$ using the standard estimator based on trajectory averages
\begin{equation}
\widehat{\vp}_{K}=\frac{1}{K+1}\sum_{k=0}^{K}\vp\left(\theta_k, p_k\right).
\label{eq: estimator vp}
\end{equation}
It can be shown in certain cases that these discrete averages converge to an ensemble average,
\[
\lim_{K\rightarrow\infty}\widehat{\vp}_{K}=\int_{\Omega}\vp(\theta, p)\ d\pi_{\varepsilon}\left(\theta, p\right)=\E_{\pi_{\varepsilon}}\!\left(\vp\right),
\]
where $\pi_{\varepsilon}$ represents the stationary density of the (biased) discrete process.
%The total error can be decomposed into error due to the bias on the invariant measure and the statistical error:
%\begin{equation}
%\widehat{\vp}_{N_{\rm iter}}-\E_{\pi}\left(\vp\right)=\Big(\E_{\pi_{\varepsilon}}\left(\vp\Big)-\E_{\pi}\left(\vp\right)\right)+\Big(\widehat{\vp}_{N_{\rm iter}}-\E_{\pi_{\varepsilon}}\left(\vp\right)\Big).
%\label{eq: error averages decomposition}
%\end{equation}
Under specific assumptions on the splitting scheme~(\cite{leimkuhler2015computation}), an expansion may be made in the time stepsize of the invariant measure of the splitting scheme which guarantees that the error in an ergodic approximation of an observable average  is bounded relative to $\varepsilon^q$, where $q$ depends on the detailed structure of the numerical method.  Thus 
\[
\E_{\pi_{\varepsilon}}\!\left(\vp\right)
=
\E_{\pi_{\rm ext}}\!\left(\vp\right) + O(\varepsilon^q).
\]
Schemes such as ``ABO'' can be shown to be first order ($q =1$), whereas BAOAB and ABOBA are second order.   Delicate cancellations imply that BAOAB can exhibit an unexpected fourth order of accuracy in the ``high friction'' limit ($\gamma\rightarrow \infty$) when the target is sampling of configurational ($\theta$-dependent) quantities.   The latter method also has remarkable features with respect configurational averages in harmonic systems and near-harmonic systems \cite{Leimkuhler2015}.

All explicit Langevin integrators are subject to stability restrictions which require that the product of stepsize and the frequency of the fastest oscillatory mode is bounded. 

\paragraph{Hybrid Monte Carlo}
	
Hybrid Monte-Carlo, also called Hamiltonian Monte Carlo, is a \glostext{MCMC} method based on the \glosname{MH} algorithm that allows one to sample directly from the target distribution $\pi$. Starting from an initial position $(\theta, p)$, the momenta are re-sampled from $\mathcal{N}\left(0,M\beta^{-1}\right)$ and a proposal $(\tilde{\theta}, \tilde{p})$ is obtained by evaluating $K_{\rm HMC}$ times the Verlet method, i.e., $B_{\varepsilon/2}A_{\varepsilon}B_{\varepsilon/2}$.  The proposal is then accepted with a probability given by the Metropolis ratio:
\begin{equation}
\min \left(1,{\rm e}^{-\beta\left(H (\tilde{\theta}, \tilde{p})- H\left(\theta, p \right)\right)}\right)\,,
\label{rejection rate hmc}
\end{equation}
where the energy $H$ is given by $H(\theta, p) = L(\theta, {\cal D}) + \frac{1}{2}p^T M^{-1}p$. In  case  the proposal is rejected, the state is reset to the starting point $(\theta, p)$. The number of steps $K_{\rm HMC}$ is often randomized. 
The parameterization of this method depends on the trade-off between the larger time stepsizes $\varepsilon$ and the number of steps $K_{\rm HMC}$, implying higher rejection rates for the proposals and smaller stepsizes leading to potentially slow exploration of the parameter space.  To our knowledge, there is currently no extension of HMC to the case of noisy gradients which retains the exact sampling property.

\paragraph{Multiple walkers: the ensemble quasi-Newton method}
We next comment on more exotic sampling schemes which can give higher sampling efficiency.   
One such scheme is the preconditioned method developed in \cite{Matthews2018} which we refer to as the \enquote{\glosname{EQN}} (\glsentrytext{glos:EQN}) method. The idea of this scheme is easily motivated by reference to a simple harmonic model problem in two space dimensions in which the stiffness (or frequency of oscillation) is very high in one direction and not in the other.  The stepsize for stable simulation using a Langevin dynamics strategy will be determined by the high frequency term meaning that the exploration rate suffers in the slowest direction. Effectively the integrated autocorrelation time in the \enquote{slow} direction  is large compared to the timestep of dynamics.

In the harmonic case it is easy to rescale the dynamical system in such a way that sampling proceeds rapidly.     In the more complicated setting of nonlinear systems in multiple dimensions, the idea that a few directions may restrict the progress of the sampling still has merit, but we no longer have direct access to the underpinning frequencies.  The idea is to determine this rescaling adaptively and dynamically during simulation.   There are several ways to do this in practice (see \cite{LeCun1991} for an early related work on using Hessian information to speed up convergence of \glostext{GD}).

\begin{figure}[h]
\begin{center}
    \includegraphics[width=4in]{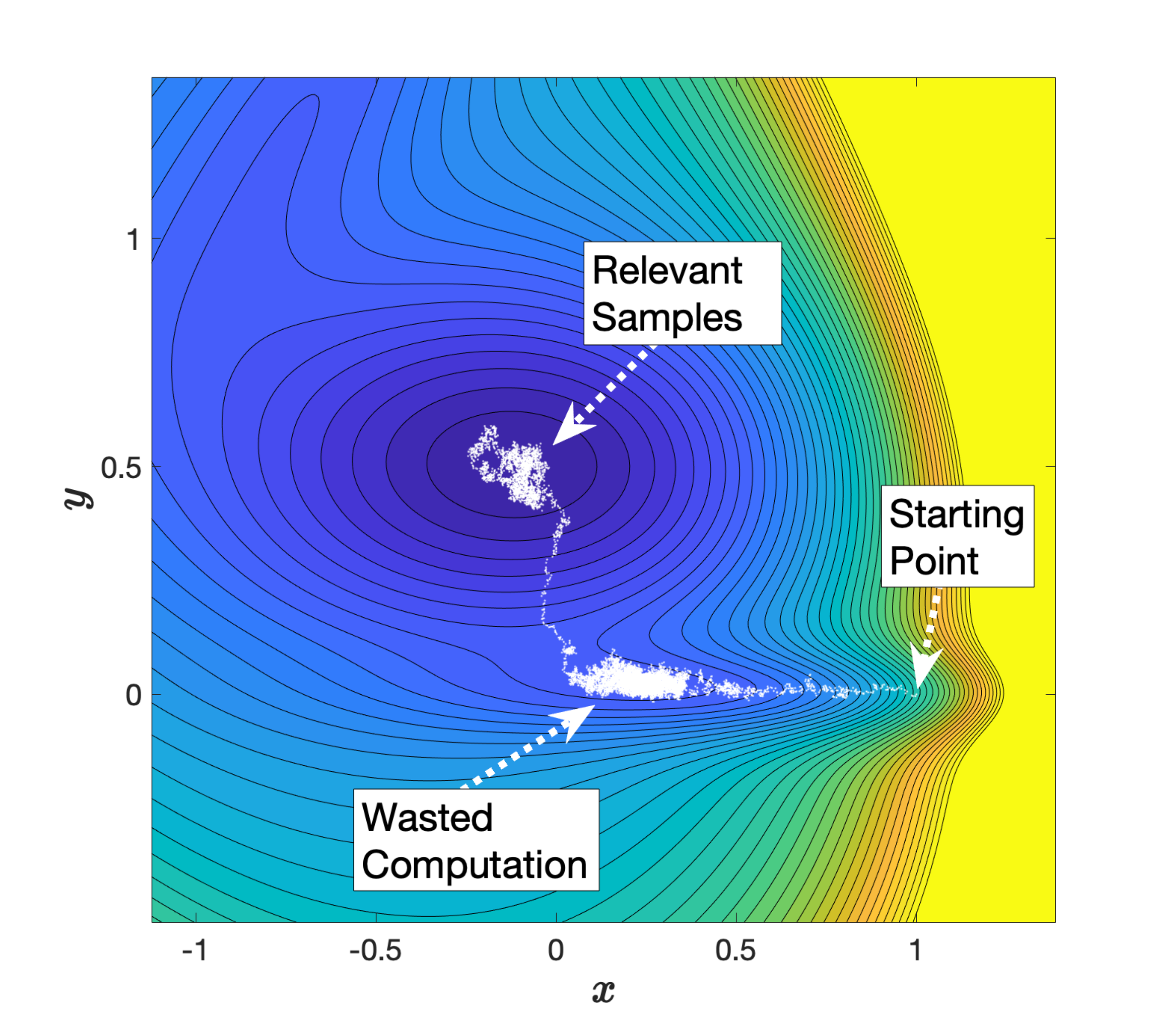}
    \caption{Slow exploration of a complex landscape using a (naive) canonical sampler.   Here Euler-Maruyama has been applied to overdamped Langevin dynamics to sample a Gaussian mixture model on a two-dimensional state space.  A poor initialization in this case leads to the extended sampling of an irrelevant intermediate basin (it contributes little to the sampling of the  overall canonical measure).  The consequence of poor scaling of the intermediate region is wasted computational effort. The goal of the enhanced sampling procedures is to accelerate the exploration in poorly scaled domains. \label{slow-sampling} }
\end{center}
\end{figure}

As a simple illustration  consider a Gausssian mixture model as shown in Figure \ref{slow-sampling}.  Here, a poor choice of initial condition has led to a naive sampling path (here generated using Euler-Maruyama and shown in white) getting stuck for a long time in a poorly scaled basin.  Eventually the sampler proceeds to the deeper (and more relevant) basin, but not before a great deal of useless computational effort has been expended.  Note that the stepsize threshold for stable integration of the SDE is inversely proportional to the frequency of the largest normal mode of oscillation in the local basin, thus it is not possible to simply 'step over' the irrelevant intermediate region by using large stepsizes.   This process mimics the behavior of many optimization schemes as well.

 In the paper of \cite{Matthews2018} modified dynamics ("Ensemble Quasi-Newton" or EQN for short) is based on construction of the covariance matrix of the walker collection.  This matrix is computed during simulation to determine the appropriate dynamical rescaling.   The implementation is somewhat involved, thus making  it an excellent demonstration of \glostext{TATi}'s versatility and robustness.   

We next briefly describe the EQN scheme; for more detail, the reader is referred to \cite{Matthews2018} (and, in particular, the \Pyth{} code referenced within it).     Suppose we have $L$ walkers (replicas).
Denote by $\theta_i$ the position vector of the $i$th walker and by $p_i$ the corresponding momentum vector, each of which is a vector in $\reel^N$, where $N$ is the number of parameters of the network. We assume that the mass matrix of the underlying system is the identity matrix, for simplicity.  The equations of motion for the $i$th walker then take the form
\begin{eqnarray}
d\theta_i & = & B_i(\theta) p_i d t,\\
d{p}_i & = & -B_i(\theta)^T\nabla L(\theta_i)dt  + {\rm div}\left ( B_i(\theta)^T \right ) d t - \gamma p_i dt + \sqrt{2\beta^{-1}\gamma} d{W}_i.
\end{eqnarray}
Here $dW_i$ has the previous interpretation of a Wiener increment and ${\rm div}$ is the tensor divergence.   The matrix $B_i$ estimates the matrix square root of the  covariance matrix  which depends on the locations of all the walkers.

Under discretization, we compute steps in configurations $\theta_{k,i}$  and momenta $p_{k,i} \in \reel^N$ at iteration step $k$ with walker index $i$.  We use a variant of the  BAOAB discretization applied to walker $i$:
\begin{subequations}
	\label{math:baoab-preconditioned}
\begin{align}
	\label{math:baoab-preconditioned-B1}
	p_{k+1/2,i} &= p_{k,i} - \frac{\varepsilon}{2} {B}_{k,i} \nabla L(\theta_{k,i}), \\
	\label{math:baoab-preconditioned-A1}
	\theta_{k+1/2, i} &= \theta_{k,i} + \frac{\varepsilon}{2} {B}_{k+1/2,i} \, p_{k+1/2,i}, \\
	\label{math:baoab-preconditioned-O}
	\hat{p}_{k+1/2,i} &=  \alpha p_{k+1/2,i} + \frac{(\alpha +1) \varepsilon}{2} \text{div} \Bigl ( B_{k+1/2,i}^T \Bigr ) + \sqrt{\frac{1-\alpha^2}{\beta}} R_{k,i}, \\
	\label{math:baoab-preconditioned-A2}
	\theta_{k+1,i} &= \theta_{k+1/2,i} + \frac{\varepsilon}{2} B_{k+1/2,i} \, \hat{p}_{k+1/2,i}, \\
	\label{math:baoab-preconditioned-B2}
	p_{k+1,i} &= \hat{p}_{k+1/2,i} - \frac{\varepsilon}{2} \nabla L(\theta_{k+1,i}).
\end{align}
\end{subequations}
As in the code of \cite{Matthews2018}, we make the calculation in \eqref{math:baoab-preconditioned-A1}  explicit by updating $B_{k,i}$  only infrequently, typically every \num{1000} steps, and we define
$B_{k,i}$ by
\begin{equation}
	\label{math:preconditoning_choice}
	B_{k,i} = \sqrt{ 
	   \unity + 
	   \eta 
	   \cov
	   \left ( \theta_{k,[i]}, \theta_{k,[i]} \right )
	   },
\end{equation}
In the above, the matrix square root is understood in the sense of a Cholesky factorization, $\eta \geq 0$ is a covariance blending constant, and $\theta_{k,[i]}$ is the set of all walker positions at timestep $k$, excluding those of walker $i$ itself.
Note that for moderate $\eta$ the choice~\eqref{math:preconditoning_choice} is always positive definite even if $L < N$. However, with \eqref{math:preconditoning_choice} the affine invariance property, an elegant feature of the original system, no longer strictly holds. Nonetheless, the simplified method has been shown to work very well in practice~\cite[\citesection~4]{Matthews2018}.

%   A number of simplifications are introduced, as in \cite{matthews2018langevin}: 
%In the context of noisy gradient, SGHMC\cite{chen2014stochastic} was proposed. 

%\todo{fh: Is the discussion of these methods that are \emph{not} contained in TATi necessary? They have been mentioned in the intro. } 
% We should discuss the stochastic versions as SGHMC\cite{chen2014stochastic}

% Riemanian Langevin ~\cite{patterson2013stochastic}

% Nogin \cite{matthews2018langevin}

% Adaptive Langevin ~\cite{shang2015covariance}

%In the case when the force is noisy~\eqref{noisy gradient}, Langevin dynamics~\eqref{eq: Langevin} the invariant                                                                                            

%\begin{equation}
  %\label{eq: Langevin}
  %\left\{
 % \begin{aligned}
%\ds    d\theta_t & =  M^{-1}p_t \, dt, \\
 %\ds   dp_t & = -\n L(\theta_t) \, dt  - \e d\tilde{W}_t - \gamma  M^{-1}p_t \, dt + \sigma \, dW_t,
  %\end{aligned}
%\right.
%\end{equation}
%the first order expansion of the invariant measure in $\e$ holds\cite{}MatthewsLeimkuhlerStoltz true:  there exists $\e^*, \varepsilon**>0$, $L > 0$ and functions $f, g$, such that, for any $0 < %\dt \leq \varepsilon^*$ and $\e \in [-\e^*, \e^*]$,
%\begin{equation}
 % \label{eq:nonequlibtrium}
  %\int_{\El} \varphi \  d \pi_{\varepsilon, \e} = \int_{\El} \varphi \ d \pi + \e \int_{\El} \varphi g \ d \pi +\dt^{\alpha} \int_{\El} \varphi f_{\alpha+1} \ d \pi + \dt^{\alpha+1} R_{\varphi,\dt}.
%\end{equation}

%\TODO	Bias, asymptotic variance, integrated auto-correlation time, KL divergence, etc

\section{Numerical Experiments}\label{sec:application}

% I have the following story in mind for the application (+ done, - todo):
% + flat potential: Asymptotic error, CLT and relation to gamma and step_width
% + quadratic potential: Arising discretization error due to non-zero gradients
% - looking at virials in QP case: what do we see here?
% + sampler properties in simple dataset: average kinetic energy and virials over all samplers
% ? sampler properties in MNIST: flat potential at high temperature, good results at low temperatures?
% + Continuing with MNIST analysis: Loss manifold to underline the large scale/high temperature, small scale/low temperature picture
% + Preliminary on EQN: Gaussian Mixture model
% - Inspecting large scale funnel: EQN for MNIST

In this section, we first look at the convergence rates and accuracy of the samplers described previously for the simplified case of a harmonic oscillator,  comparing them with analytically known rates.  Next we turn to a very simple clustering problem to further explore the error (bias) introduced by SDE schemes; we show that a particular discretization (``BAOAB'') of underdamped Langevin dynamics offers extraordinarily high accuracy compared to several alternative, mimicing observations about this scheme in the molecular dynamics setting.
Subsequently, we use the MNIST dataset on a single-layer perceptron to illustrate an enhanced loss landscape visualization technique that obtains its projection directions from sampling trajectories.
Finally, we present results on the \glostext{EQN}  method described in Section~\ref{sec:sampling} for a Gaussian model and for a single layer perceptron applied to the MNIST dataset.

\subsection{Sampler Properties and Error Analysis}\label{sec:application-samplers}
In this error analysis, we explore a state $\{\theta,p\}$ with position $\theta$ and momentum $p$ of a single degree of freedom. The  kinetic energy is defined as $\varphi(\theta, p) = \tfrac 1 2 p^T p$, where we have set the mass matrix to unity. Its asymptotic value in the canonical distribution is given by the number of parameters $N$ and the (inverse) temperature $\beta$ as $\tfrac{1}{2\beta}$   \cite[\citesection~6.1.5]{Leimkuhler2012}.
The virial is defined as $\varphi(\theta, p) = \sum^N_i \theta_i \cdot \nabla_i L(\theta) = \theta \cdot \nabla L(\theta)$, where $\nabla_i L(\theta)$ is the derivative of the loss function with respect to the parameter $\theta_i$. Its asymptotic value is two times that of the kinetic energy, as a consequence of the virial theorem. For this theorem to hold we need a potential that is unbounded from above and grows sufficiently rapidly at infinity, see appendix \ref{sec:virial_theorem-softmax_cross_entropy} for a derivation.

Averages of the kinetic energy or of the virial are examples of time integrals of \eqref{eq: ergodic averages markov chain} that can be used to assess the accuracy of the chosen dynamics and discretization, since their asymptotic values are known.  As we are interested in discretized integrals over finite number of time steps, which are accessible in numerical computations in the absence of analytical solutions,  we look at estimators~\eqref{eq: estimator vp} of the following form
\[
	\hat{\varphi}_N = \frac 1 N \sum^{N-1}_{n=0} \varphi(q_n, p_n).
\]
The total error with respect to the expected value $E_\mu(\varphi)$ for the invariant probability measure $\mu$ decomposes as
\begin{equation}
	\label{eq:total_error_estimator}
	E\left ( | \hat{\varphi}_N - E_\mu (\varphi) |^2 \right ) = \left ( E(\hat{\varphi}_N)- E_\mu(\varphi) \right )^2 + E\left ( | \hat{\varphi}_N - E(\hat{\varphi}_N) |^2\right ),
\end{equation}
i.e. it is a combination  of the \emph{discretization error},  from the finite step size when discretizing the dynamics~\eqref{eq: Langevin}, and the \emph{sampling error} that results from the inability to sample over an infinite time or to generate an infinite number of steps. The first source is also sometimes referred to as the  perfect sampling bias, emphasizing its presence even in the limit of infinitely many samples.

In the extreme case of a vanishing gradient, only the sampling error is present. This error will decay so that its variance is proportional to $1/N$.
A more interesting case is that of a \enquote{harmonic potential} $I(\theta) = a \theta^2/2$, in which case the truncation error is nontrivial.

\subsubsection{Truncation Error for Harmonic Potential}

%Before, we have investigated the behavior of the sampler as a random walker on a flat potential. We have seen how the error of the average kinetic energy to its asymptotic value changes with the friction constant $\gamma$ and the step size $\varepsilon$.
%Let us now introduce a harmonic potential in the form $l(\theta) = a \cdot \theta^2$, where $a$ acts as a scale to the potential.

We have contributions to the total error~\eqref{eq:total_error_estimator} from both the sampling error and the discretization error. We will see that the latter may dominate when the scale $a$ of the potential and therefore the average gradients are sufficiently large.  We concentrate here on the empirical results; refer to \citet[\citesection~7.4]{Leimkuhler2012} for a general discussion of harmonic problems in the context of Langevin dynamics.    We use a single-layer perceptron with a single input node and a single output node with linear activation and zero bias, i.\,e.~$f_{\theta}(x_i) = \theta_1 x_i$.  We use the mean squared loss function.
Then, such a harmonic potential can be easily introduced by a dataset with the square root of the prefactor as its single feature and a zero label, i.\,e.~the only dataset item is $(X,Y)=(\sqrt{a},0)$.
We use three different factors $a \in \{0.01, 1, 4\}$. 
Moreover, we use the following sets of parameters for $\gamma \in \{0.01, 0.1, 1, 10, 100 \}$, $\beta = 10$, and $\varepsilon \in \{ 2^{-2}, 2^{-3}, 2^{-4}, 2^{-5}, 2^{-6} \}$. Here, we employ BAOAB as the sampler, again with $N=10^6$ sampling steps.

% See quadratic potential, "15_Asymptotic_Kinetic_Error"

\begin{figure}[htbp]
	\def\myestimated{kinetic_energy}%
	\subfigure[Small prefactor $a=0.01$]{
		\label{fig:quadratic_potential-small_prefactor}
		\resizebox{0.49\textwidth}{!}{
			\def\myfactor{0_1}
			\def\myfactordot{0.1}
			\tikzsetnextfilename{step_width-error_to_asymptotic_value-slope_fit-baoab-small_prefactor}
			\ifdefined\myexperimentname\empty\else\def\myexperimentname{Quadratic potential}\fi%
\ifdefined\mypath\empty\else\def\mypath{\myexperimentspath/15_Asymptotic_Kinetic_Error/data/Quadratic_potential}\fi%

%special column names
\ifdefined\myxcolumnname\empty\else\def\myxcolumnname{step_width}\fi%
\ifdefined\myycolumnname\empty\else\def\myycolumnname{final_mean_value}\fi%
\ifdefined\myestimated\empty\else\def\myestimated{kinetic_energy}\fi%
%parameters
\ifdefined\myfactor\empty\else\def\myfactor{2_}\fi% %0_1,1_,2_
\ifdefined\myfactordot\empty\else\def\myfactordot{2.}\fi% %0.1,1.,2.\ifdefined\mygamma\empty\else\def\mygamma{0.01}\fi% %0.01, 0.1, 1, 10, 100

\pgfplotstableset{
	create on use/step_width/.style={create col/copy column},
	create on use/kinetic_energy_error/.style={create col/expr={abs(0.05-\thisrow{avg_\myycolumnname})/0.05}},
}

\begin{tikzpicture}[scale=1.2]
\begin{loglogaxis}[
	width=7cm,height=6cm,
	%title style={align=left},
	%title = {\emph{\myexperimentname}}: \expandafter\spaceus\expandafter{\myestimated} with BAOAB,%:\\ dashed - true final value, solid - estimated final value},
	xlabel = {Step width},
	ylabel = {Average kinetic energy error},
	xtick={0.015625,0.03125,0.0625,0.125,0.25},
	log ticks with fixed point,
	legend pos=outer north east,
	grid=major,
%	each nth point=1,
	cycle multi list={
		{mark=o,blue},{mark=square,red},{mark=otimes, brown},{mark=star,black},{mark=diamond,green}\nextlist
		dashed,solid
	},
]
\addlegendimage{empty legend}
\addlegendentry{$\gamma$}
\foreach \mygamma in {0.01,0.1,1,10,100}{%
% \addplot+[error bars/.cd, y dir=both, y explicit] table[x=\myxcolumnname,y expr={abs(0.05-\thisrow{avg_final_mean_value})/0.05}, y error expr={sqrt(\thisrow{var_final_mean_value}-\thisrow{avg_final_mean_value}^2)/0.05}, col sep=comma]{\mypath/estimated-factor_\myfactordot-gamma_\mygamma.csv};
% \expandafter\addlegendentry\expandafter{\mygamma}

\pgfplotstableread[columns={step_width,average_\mycolumnname}, col sep=comma]{\mypath/estimated_\myestimated-factor_\myfactordot-gamma_\mygamma.csv}\mybaoabtable

\addplot+[thick,dotted] table[x=\myxcolumnname, y={create col/linear regression={x=\myxcolumnname,y=kinetic_energy_error, variance list={100,10,1,1,1}}}]\mybaoabtable;
\xdef\slopeglaone{\pgfplotstableregressiona}
\addlegendentryexpanded{slope: \pgfmathprintnumber[fixed,precision=2]{\slopeglaone}}

\addplot+[thick,error bars/.cd, y dir=both, y explicit] table[x=\myxcolumnname,y expr={abs(0.05-\thisrow{avg_\myycolumnname})/0.05}, y error expr={sqrt(\thisrow{var_\myycolumnname}-\thisrow{avg_\myycolumnname}^2)/0.05}, col sep=comma]{\mypath/estimated_\myestimated-factor_\myfactordot-gamma_\mygamma.csv};
\expandafter\addlegendentry\expandafter{\mygamma}
}%
\end{loglogaxis}
\end{tikzpicture}
		}
	}
	\subfigure[Large prefactor $a=4$]{
		\label{fig:quadratic_potential-large_prefactor}
		\resizebox{0.49\textwidth}{!}{
			\def\myfactor{2_}
			\def\myfactordot{2.}
			\tikzsetnextfilename{step_width-error_to_asymptotic_value-slope_fit-baoab-large_prefactor}
			\ifdefined\myexperimentname\empty\else\def\myexperimentname{Quadratic potential}\fi%
\ifdefined\mypath\empty\else\def\mypath{\myexperimentspath/15_Asymptotic_Kinetic_Error/data/Quadratic_potential}\fi%

%special column names
\ifdefined\myxcolumnname\empty\else\def\myxcolumnname{step_width}\fi%
\ifdefined\myycolumnname\empty\else\def\myycolumnname{final_mean_value}\fi%
\ifdefined\myestimated\empty\else\def\myestimated{kinetic_energy}\fi%
%parameters
\ifdefined\myfactor\empty\else\def\myfactor{2_}\fi% %0_1,1_,2_
\ifdefined\myfactordot\empty\else\def\myfactordot{2.}\fi% %0.1,1.,2.\ifdefined\mygamma\empty\else\def\mygamma{0.01}\fi% %0.01, 0.1, 1, 10, 100

\pgfplotstableset{
	create on use/step_width/.style={create col/copy column},
	create on use/kinetic_energy_error/.style={create col/expr={abs(0.05-\thisrow{avg_\myycolumnname})/0.05}},
}

\begin{tikzpicture}[scale=1.2]
\begin{loglogaxis}[
	width=7cm,height=6cm,
	%title style={align=left},
	%title = {\emph{\myexperimentname}}: \expandafter\spaceus\expandafter{\myestimated} with BAOAB,%:\\ dashed - true final value, solid - estimated final value},
	xlabel = {Step width},
	ylabel = {Average kinetic energy error},
	xtick={0.015625,0.03125,0.0625,0.125,0.25},
	log ticks with fixed point,
	legend pos=outer north east,
	grid=major,
%	each nth point=1,
	cycle multi list={
		{mark=o,blue},{mark=square,red},{mark=otimes, brown},{mark=star,black},{mark=diamond,green}\nextlist
		dashed,solid
	},
]
\addlegendimage{empty legend}
\addlegendentry{$\gamma$}
\foreach \mygamma in {0.01,0.1,1,10,100}{%
% \addplot+[error bars/.cd, y dir=both, y explicit] table[x=\myxcolumnname,y expr={abs(0.05-\thisrow{avg_final_mean_value})/0.05}, y error expr={sqrt(\thisrow{var_final_mean_value}-\thisrow{avg_final_mean_value}^2)/0.05}, col sep=comma]{\mypath/estimated-factor_\myfactordot-gamma_\mygamma.csv};
% \expandafter\addlegendentry\expandafter{\mygamma}

\pgfplotstableread[columns={step_width,average_\mycolumnname}, col sep=comma]{\mypath/estimated_\myestimated-factor_\myfactordot-gamma_\mygamma.csv}\mybaoabtable

\addplot+[thick,dotted] table[x=\myxcolumnname, y={create col/linear regression={x=\myxcolumnname,y=kinetic_energy_error, variance list={100,10,1,1,1}}}]\mybaoabtable;
\xdef\slopeglaone{\pgfplotstableregressiona}
\addlegendentryexpanded{slope: \pgfmathprintnumber[fixed,precision=2]{\slopeglaone}}

\addplot+[thick,error bars/.cd, y dir=both, y explicit] table[x=\myxcolumnname,y expr={abs(0.05-\thisrow{avg_\myycolumnname})/0.05}, y error expr={sqrt(\thisrow{var_\myycolumnname}-\thisrow{avg_\myycolumnname}^2)/0.05}, col sep=comma]{\mypath/estimated_\myestimated-factor_\myfactordot-gamma_\mygamma.csv};
\expandafter\addlegendentry\expandafter{\mygamma}
}%
\end{loglogaxis}
\end{tikzpicture}
		}
	}
	\caption{Highlighting the sampling error with respect to the step size by showing the relative error of the average kinetic energy with respect to its asymptotic value for the \enquote{quadratic potential} case. \subref{fig:quadratic_potential-small_prefactor}{} If the gradients are small, the system performs a pure random walk and results are independent of the step size. \subref{fig:quadratic_potential-large_prefactor}{} Given large enough gradients relative to the temperature the convergence order of the integration method with the step size becomes apparent, here second order for BAOAB in the momenta.}
	\label{fig:quadratic_potential}
\end{figure}
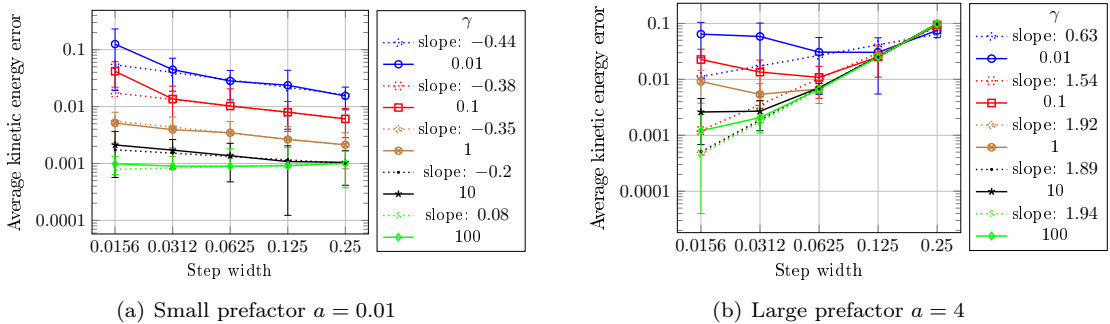

% Potential gives rise to discretization error

Examining Figure~\ref{fig:quadratic_potential-small_prefactor} where we use a small prefactor $a$, we notice that the error decreases with increasing step size because $\alpha$ becomes smaller and therefore we obtain more random walk-like behavior. However, it does not entirely depend on $\alpha$, but also to some extent on the step size $\varepsilon$. For the highest value of $\gamma=100$ we again have a flat line due to the lower bound enforced by the \glostext{CLT}.  

In Figure~\ref{fig:quadratic_potential-large_prefactor} with a large prefactor $a$ for large step sizes all of the curves coincide regardless of $\gamma$ and the behavior has reversed: now the error becomes smaller for smaller step sizes.  Naturally, the reason for this change is the discretization error that arises because of substantial non-zero gradients, and that the error depends on the step size $\varepsilon$. Measuring the slope in the domain where all curves overlap for $a=4$, we obtain values of up to $2$, i.\,e.~second order convergence in the discretization error as expected from BAOAB.  In place of the prefactor $a$ we could also have varied the inverse temperature $\beta$ to the same effect that only depends on the scale of the noise relative to the scale of the gradients.

%Hence, we have seen that there are two error sources for this average kinetic energy: On the one hand it is the finite trajectory length and on the other hand it is the discretization and its finite step size. 

Using a higher-order sampler allows for a smaller error at a given step size or to use larger step sizes (and therefore sample more space) for a given error threshold. Note that the step size is bounded from above by a stability threshold.   The benefit of the trade-off between accuracy and computational effort is limited however, and it is likely that very high order schemes (beyond the second order splittings discussed here) are not efficacious in TATi, in keeping with previous studies in molecular dynamics (\cite{Leimkuhler2015}).

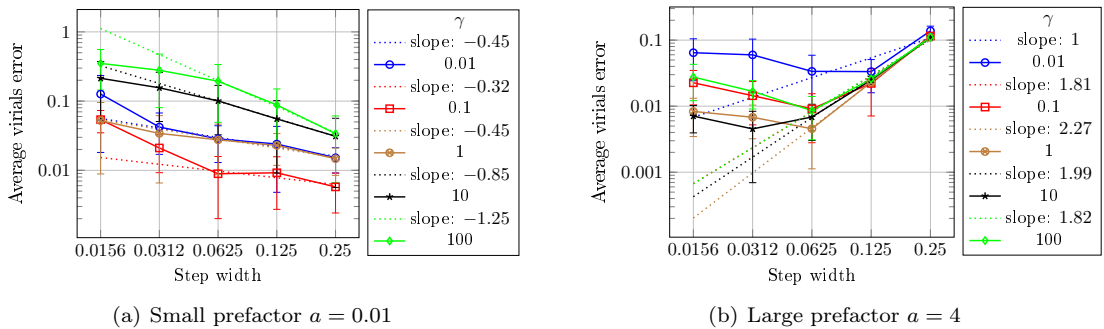
\begin{figure}[htbp]
	\def\myestimated{virial}%
	\subfigure[Small prefactor $a=0.01$]{
		\label{fig:quadratic_potential_gla2-small_prefactor}
		\resizebox{0.49\textwidth}{!}{
			\def\myfactor{0_1}
			\def\myfactordot{0.1}
			\tikzsetnextfilename{step_width-error_to_asymptotic_value-slope_fit-gla2-small_prefactor}
			\ifdefined\myexperimentname\empty\else\def\myexperimentname{Quadratic potential - GLA2}\fi%
\ifdefined\mypath\empty\else\def\mypath{\myexperimentspath/15_Asymptotic_Kinetic_Error/data/Quadratic_potential_-_GLA2}\fi%

%special column names
\ifdefined\myxcolumnname\empty\else\def\myxcolumnname{step_width}\fi%
\ifdefined\myycolumnname\empty\else\def\myycolumnname{final_mean_value}\fi%
\ifdefined\myestimated\empty\else\def\myestimated{virial}\fi%
%parameters
\ifdefined\myfactor\empty\else\def\myfactor{2_}\fi% %0_1,1_,2_
\ifdefined\myfactordot\empty\else\def\myfactordot{2.}\fi% %0.1,1.,2.\ifdefined\mygamma\empty\else\def\mygamma{0.01}\fi% %0.01, 0.1, 1, 10, 100

\pgfplotstableset{
	create on use/step_width/.style={create col/copy column},
	create on use/kinetic_energy_error/.style={create col/expr={abs(0.05-\thisrow{avg_\myycolumnname})/0.05}},
}

\begin{tikzpicture}[]
\begin{loglogaxis}[
	width=7cm,height=6cm,
	%title style={align=left},
	%title = {\emph{\myexperimentname}}: \expandafter\spaceus\expandafter{\myestimated},%:\\ dashed - true final value, solid - estimated final value},
	xlabel = {Step width},
	ylabel = {Average virials error},
	xtick={0.015625,0.03125,0.0625,0.125,0.25},
	legend pos=outer north east,
	log ticks with fixed point,
	grid=major,
%	each nth point=1,
	cycle multi list={
		{mark=o,blue},{mark=square,red},{mark=otimes, brown},{mark=star,black},{mark=diamond,green}\nextlist
		dashed,solid
	},
]
\addlegendimage{empty legend}
\addlegendentry{$\gamma$}
\foreach \mygamma in {0.01,0.1,1,10,100}{%
% \addplot+[error bars/.cd, y dir=both, y explicit] table[x=\myxcolumnname,y expr={abs(0.05-\thisrow{avg_final_mean_value})/0.05}, y error expr={sqrt(\thisrow{var_final_mean_value}-\thisrow{avg_final_mean_value}^2)/0.05}, col sep=comma]{\mypath/estimated-factor_\myfactordot-gamma_\mygamma.csv};
% \expandafter\addlegendentry\expandafter{\mygamma}

\pgfplotstableread[columns={step_width,average_\mycolumnname}, col sep=comma]{\mypath/estimated_\myestimated-factor_\myfactordot-gamma_\mygamma.csv}\mybaoabtable

\addplot+[thick,dotted,no marks] table[x=\myxcolumnname, y={create col/linear regression={x=\myxcolumnname,y=kinetic_energy_error, variance list={100,10,1,1,1}}}]\mybaoabtable;
\xdef\slopeglaone{\pgfplotstableregressiona}
\addlegendentryexpanded{slope: \pgfmathprintnumber[fixed,precision=2]{\slopeglaone}}

\addplot+[thick,error bars/.cd, y dir=both, y explicit] table[x=\myxcolumnname,y expr={abs(0.05-\thisrow{avg_\myycolumnname})/0.05}, y error expr={sqrt(\thisrow{var_\myycolumnname}-\thisrow{avg_\myycolumnname}^2)/0.05}, col sep=comma]{\mypath/estimated_\myestimated-factor_\myfactordot-gamma_\mygamma.csv};
\expandafter\addlegendentry\expandafter{\mygamma}
}%
\end{loglogaxis}
\end{tikzpicture}
		}
	}
	\subfigure[Large prefactor $a=4$]{
		\label{fig:quadratic_potential_gla2-large_prefactor}
		\resizebox{0.49\textwidth}{!}{
			\def\myfactor{2_}
			\def\myfactordot{2.}
			\tikzsetnextfilename{step_width-error_to_asymptotic_value-slope_fit-gla2-large_prefactor}
			\ifdefined\myexperimentname\empty\else\def\myexperimentname{Quadratic potential - GLA2}\fi%
\ifdefined\mypath\empty\else\def\mypath{\myexperimentspath/15_Asymptotic_Kinetic_Error/data/Quadratic_potential_-_GLA2}\fi%

%special column names
\ifdefined\myxcolumnname\empty\else\def\myxcolumnname{step_width}\fi%
\ifdefined\myycolumnname\empty\else\def\myycolumnname{final_mean_value}\fi%
\ifdefined\myestimated\empty\else\def\myestimated{virial}\fi%
%parameters
\ifdefined\myfactor\empty\else\def\myfactor{2_}\fi% %0_1,1_,2_
\ifdefined\myfactordot\empty\else\def\myfactordot{2.}\fi% %0.1,1.,2.\ifdefined\mygamma\empty\else\def\mygamma{0.01}\fi% %0.01, 0.1, 1, 10, 100

\pgfplotstableset{
	create on use/step_width/.style={create col/copy column},
	create on use/kinetic_energy_error/.style={create col/expr={abs(0.05-\thisrow{avg_\myycolumnname})/0.05}},
}

\begin{tikzpicture}[]
\begin{loglogaxis}[
	width=7cm,height=6cm,
	%title style={align=left},
	%title = {\emph{\myexperimentname}}: \expandafter\spaceus\expandafter{\myestimated},%:\\ dashed - true final value, solid - estimated final value},
	xlabel = {Step width},
	ylabel = {Average virials error},
	xtick={0.015625,0.03125,0.0625,0.125,0.25},
	legend pos=outer north east,
	log ticks with fixed point,
	grid=major,
%	each nth point=1,
	cycle multi list={
		{mark=o,blue},{mark=square,red},{mark=otimes, brown},{mark=star,black},{mark=diamond,green}\nextlist
		dashed,solid
	},
]
\addlegendimage{empty legend}
\addlegendentry{$\gamma$}
\foreach \mygamma in {0.01,0.1,1,10,100}{%
% \addplot+[error bars/.cd, y dir=both, y explicit] table[x=\myxcolumnname,y expr={abs(0.05-\thisrow{avg_final_mean_value})/0.05}, y error expr={sqrt(\thisrow{var_final_mean_value}-\thisrow{avg_final_mean_value}^2)/0.05}, col sep=comma]{\mypath/estimated-factor_\myfactordot-gamma_\mygamma.csv};
% \expandafter\addlegendentry\expandafter{\mygamma}

\pgfplotstableread[columns={step_width,average_\mycolumnname}, col sep=comma]{\mypath/estimated_\myestimated-factor_\myfactordot-gamma_\mygamma.csv}\mybaoabtable

\addplot+[thick,dotted,no marks] table[x=\myxcolumnname, y={create col/linear regression={x=\myxcolumnname,y=kinetic_energy_error, variance list={100,10,1,1,1}}}]\mybaoabtable;
\xdef\slopeglaone{\pgfplotstableregressiona}
\addlegendentryexpanded{slope: \pgfmathprintnumber[fixed,precision=2]{\slopeglaone}}

\addplot+[thick,error bars/.cd, y dir=both, y explicit] table[x=\myxcolumnname,y expr={abs(0.05-\thisrow{avg_\myycolumnname})/0.05}, y error expr={sqrt(\thisrow{var_\myycolumnname}-\thisrow{avg_\myycolumnname}^2)/0.05}, col sep=comma]{\mypath/estimated_\myestimated-factor_\myfactordot-gamma_\mygamma.csv};
\expandafter\addlegendentry\expandafter{\mygamma}
}%
\end{loglogaxis}
\end{tikzpicture}
		}
	}
	\caption{Relative error of the average virial with respect to its asymptotic value for the \enquote{quadratic potential} case using \protect\glostext{GLA}2 sampler. Again, we see second order convergence, here in the positions, given the gradients are large enough relative to the noise, see also Figure~\ref{fig:quadratic_potential}.}
	\label{fig:quadratic_potential_gla2}
\end{figure}

For comparison, we also look at the average virials in Figure~\ref{fig:quadratic_potential_gla2}. However, they depend only on the position and \emph{not} on momentum.   Note that the BAOAB scheme has nil perfect sampling bias for purely configurational quantities such as virials.
Here, we are instead using the \glostext{GLA} 2nd order sampler but keep all the other aspects of the method unchanged.

% See quadratic potential - GLA2, "15_Asymptotic_Kinetic_Error"

We obtain the same qualitative picture for the average virials sampling with \glostext{GLA}2 as we got with the average kinetic energy sampling with BAOAB. Again, for a large enough prefactor the discretization error dominates. 
Inspecting the slopes in the doubly logarithmic plots, we find values around 2 that peak for $\gamma=1$. 
%Note that BAOAB did not show such a peak. This is not surprising as BAOAB is expected to even show fourth order convergence in the high-friction limit, see \cite{Leimkuhler2012}.

\subsubsection{Langevin Sampler Performance in a Two-Cluster Classification Problem}

In the following we will be inspecting the average virial, obtained over sampled, finite trajectories in order to assess the accuracies of positions obtained from various samplers. We will be investigating the following samplers: \glosfirst{SGLD}, \glosfirst{GLA} 1st and 2nd order, and BAOAB. %Note that \glosfirst{SGLD} corresponds to a (first-order) Euler-Maruyama discretization of over-damped Langevin dynamics.

Very long runs are needed to bring forth the different convergence order of the discretization error because of the involved statistical errors. Therefore, we still use a very simple dataset: it consists of 500 points drawn from two Gaussians in two dimensions, one centered at $[2,2]$ with label $1$, the other centered at $[-2,-2]$ with label $-1$. The points are additionally perturbed by 0.1 relative noise. 

We use a single layer perceptron with two input nodes, a single output node with linear activation and mean squared loss. Therefore, the network has $N=3$ parameters in total, two weight degrees and one bias degree.

The parameters are first equilibrated for \num{2000} steps with a learning rate of 0.03 with \glosfirst{GD}. Next, we perform sampling runs from the resulting position for $10^6$ steps at various step sizes $\varepsilon$. In the case of a sampler based on Langevin dynamics, we use a friction constant $\gamma=10$ and an inverse temperature $\beta=10$. Note that, because we start at an equilibrated position with zero gradients and because the potential function is squared, positions are only rescaled when using other temperatures if the same random number sequence is used. We use 100 different seeds and average the last average value per trajectory over all seeds.

\begin{figure}[htbp]
	\def\myexperimentspath{/nosave/heber/Experiments/DataDrivenSampler}
%	\subfigure[Kinetic energy]{
%		\resizebox{0.45\textwidth}{!}{
%			\tikzsetnextfilename{order_of_decay-slope_fit-all_samplers-kinetic}
%			\input{./figures/10_OrderOfDecay_TwoClusters/order_of_decay-slope_fit-all_samplers-kinetic}
%		}
%	}
%	\subfigure[Virial]{
	\centering
		\resizebox{0.9\textwidth}{!}{
			\tikzsetnextfilename{order_of_decay-slope_fit-all_samplers-virials}
			\ifdefined\mypath\empty\else\def\mypath{\myexperimentspath/../DataDrivenSampler/10_OrderOfDecay_TwoClusters/data/}\fi%

% parameter
\ifdefined\mysamplername\empty\else\def\mysampler{StochasticGradientLangevinDynamics}\fi%
\ifdefined\mysamplername\empty\else\def\mysamplername{SGLD}\fi%
\ifdefined\mybeta\empty\else\def\mybeta{10}\fi%
\ifdefined\mygamma\empty\else\def\mygamma{10}\fi%

\def\mydof{3}%
\def\myasymptoticloss{0.2349601251}%

\ifdefined\mycolumn\empty\else\def\mycolumn{virials}\fi%

% samplers

\begin{tikzpicture}[]

% SGLD
\pgfplotstableread[columns={step_width,average_\mycolumnname}, col sep=comma]{\mypath/Order_of_decay_-_SGLD/average_degrees_of_freedom-friction_constant_\mygamma-inverse_temperature_\mybeta-sampler_StochasticGradientLangevinDynamics.csv}\mysgldtable
% GLA1
\pgfplotstableread[columns={step_width,average_\mycolumnname}, col sep=comma]{\mypath/Order_of_decay_-_GLA1/average_degrees_of_freedom-friction_constant_\mygamma-inverse_temperature_\mybeta-sampler_GeometricLangevinAlgorithm_1stOrder.csv}\myglaonetable
% GLA2
\pgfplotstableread[columns={step_width,average_\mycolumnname}, col sep=comma]{\mypath/Order_of_decay_-_GLA2/average_degrees_of_freedom-friction_constant_\mygamma-inverse_temperature_\mybeta-sampler_GeometricLangevinAlgorithm_2ndOrder.csv}\myglatwotable
% BAOAB
\pgfplotstableread[columns={step_width,average_\mycolumnname}, col sep=comma]{\mypath/Order_of_decay_-_BAOAB/average_degrees_of_freedom-friction_constant_\mygamma-inverse_temperature_\mybeta-sampler_BAOAB.csv}\mybaoabtable

\pgfplotstablegetrowsof\mysgldtable%

\pgfplotstableset{
	create on use/virials_error/.style={create col/expr={abs(\mydof/2-\thisrow{avg_average_virials}*\mybeta)}},
	create on use/kinetic_energy_error/.style={create col/expr={abs(\mydof/2-\thisrow{avg_average_kinetic_energy}*\mybeta)}},
	create on use/potential_error/.style={create col/expr={abs(\myasymptoticloss-\thisrow{avg_ensemble_average_loss})}},
}

% \pgfplotstablecreatecol[
% 	create on use/error/.style={create col/expr={abs(\mydof/2-\thisrow{avg_average_\mycolumnname}*1e1)}},
% ]{error}{\myglafirsttable}
%\pgfplotstablecreatecol[linear regression]{linear_regression}{\myglafirsttable}
%\xdef\slopeGLAfirst{\pgfplotstableregressiona}

% \pgfplotstablecreatecol{error}{\myglasecondtable}
%\pgfplotstablecreatecol[linear regression]{error}{\myglasecondtable}
%\xdef\slopeGLAsecond{\pgfplotstableregressiona}

% \node (a)  {\pgfplotstabletypeset[columns={step_width}]\myglafirsttable};
\begin{loglogaxis}[
	width=12cm,height=10.5cm,
	xlabel = {Step size},
	ylabel = {$\Bigl | \tfrac {N_{\text{DOF}}} 2 - E^{\text{avg}}_{\text{vir}}
\cdot \beta \Bigr |$},
	grid=major,
	legend pos=outer north east,
	legend columns=2,
]
% SGLD/Euler-Maruyama
\addplot+[thick,red,solid, mark=o] table[x expr={sqrt(\thisrow{step_width})},y=\mycolumn_error]\mysgldtable;
\addlegendentry{}

\addplot+[thick,red, dashed,no marks] table[x expr={sqrt(\thisrow{step_width})}, y={create col/linear regression={=step_width,y=\mycolumn_error, variance list={1,1,1,1,10}}}]\mysgldtable;
\xdef\slopesgld{\pgfplotstableregressiona}
%\addlegendentry{slope SGLD: \pgfmathprintnumber{\slopesgld}}
\addlegendentry{SGLD}

% GLA1
\addplot+[thick,blue,solid, mark=o] table[x=step_width,y=\mycolumn_error]\myglaonetable;
\addlegendentry{}

\addplot+[thick,blue, dashed,no marks] table[x=step_width, y={create col/linear regression={x=step_width,y=\mycolumn_error, variance list={1,1,1,1,10}}}]\myglaonetable;
\xdef\slopeglaone{\pgfplotstableregressiona}
%\addlegendentry{slope GLA1: \pgfmathprintnumber{\slopeglaone}}
\addlegendentry{GLA1}

% GLA2
\addplot+[thick,brown,solid, mark=o] table[x=step_width,y=\mycolumn_error]\myglatwotable;
\addlegendentry{}

\addplot+[thick,brown, dashed,no marks] table[x=step_width, y={create col/linear regression={x=step_width,y=\mycolumn_error, variance list={1,1,1,1,10}}}]\myglatwotable;
\xdef\slopeglatwo{\pgfplotstableregressiona}
%\addlegendentry{slope GLA2: \pgfmathprintnumber{\slopeglatwo}}
\addlegendentry{GLA2}

% BAOAB
\addplot+[thick,black,solid, mark=o] table[x=step_width,y=\mycolumn_error]\mybaoabtable;
\addlegendentry{}

\addplot+[thick,black, dashed,no marks] table[x=step_width, y={create col/linear regression={x=step_width,y=\mycolumn_error, variance list={1,1,1,1,10}}}]\mybaoabtable;
\xdef\slopebaoab{\pgfplotstableregressiona}
%\addlegendentry{slope BAOAB: \pgfmathprintnumber{\slopebaoab}}
\addlegendentry{BAOAB}
\end{loglogaxis}
\end{tikzpicture}
		}
%	}
	\caption{Order of convergence for the discretization error for four samplers of the average average virial for the simple two clusters dataset. Dynamics become unstable if the largest step sizes are increased by a factor of two. \protect\glostext{GLA}{}1 has first order (slope 0.86). \protect\glostext{SGLD} (slope 2.18) and \protect\glostext{GLA}{}2 (slope 2.01) have second order. BAOAB's accuracy on the position marginal is so good that its second to fourth order convergence cannot be seen against the lower bound of the \protect\glostext{CLT}.}
	\label{fig:order_of_decay}
\end{figure}
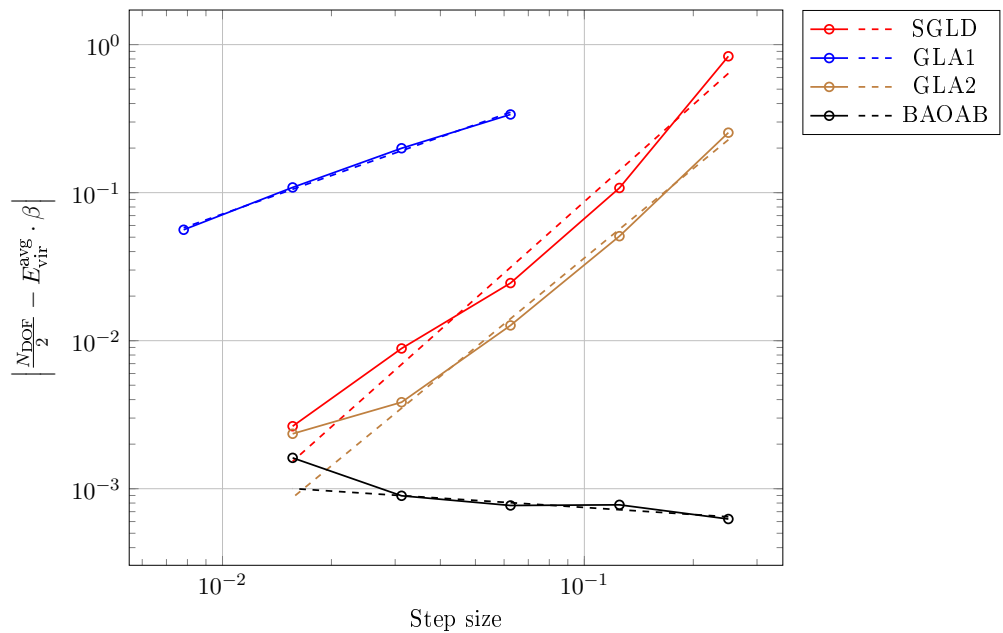

In Figure~\ref{fig:order_of_decay} we give the absolute error between the average virial and its asymptotic value scaled by the inverse temperature $\beta$ per sampler for each step size employed. Note again that the virial depends on the positions and the gradient. % while the average kinetic energy depends only on the far less important momenta.
Moreover, we remark that increasing the largest step size shown for each sampler individually by a factor of two would cause the dynamics to become unstable. Naturally, the exact threshold depends on the magnitude of the gradients and therefore on the dataset.

The slopes have been obtained from least squares regression fits where the data point to the smallest step size have been weighted by $\tfrac{1}{10}$.

We make the following observations: \glostext{GLA}{2} has second order convergence in the average virial, \glostext{GLA}{1} has first-order convergence. BAOAB shows such great accuracy at this finite trajectory length that its second to fourth order convergence does not show as it reaches the \glostext{CLT} limit. These results are in absolute agreement with results from an analysis on harmonic problems, see \cite[\citesection~7.4.1]{Leimkuhler2012}.
Note that {SGLD} exhibits second-order convergence in the virial in this example; actually this is an artifact of the way the data has been graphed: the SGLD 'step size' can be viewed as the square root of the step size $\sqrt{\varepsilon}$ of the other samplers, see \cite[\citepage~36]{Leimkuhler2012}, i.e. it is really a first order scheme when expressed in the standard way. 

%Moreover, \glostext{GLA}{2}, that can also be seen as OBABO with unified O steps (apart from the different variance), samples the momenta so accurately that we would need much longer sampling runs to see any convergent behavior in the average kinetic energy. The same holds for the virial for BAOAB, i.\,e.~it sees a high very high accuracy with the positions. Both of which is known, too.

From these results the sampling method of choice seems to be BAOAB which has superior accuracy in the positions and general second order convergence at little computational overhead and extra memory requirement, compared to \glostext{SGLD}. Note that \glostext{SGLD} is based on Brownian dynamics, thus lacking momenta, and therefore is expected to be less efficient in exploration for multimodal problems than BAOAB or other Langevin dynamics samplers~\cite{Leimkuhler2015}.

A similar study on the MNIST dataset is impeded by its wide-spread covariance eigenvalue spectrum, shown in the next section. There we encounter both large and small gradients and therefore have a mixture of the \enquote{flat potential} and \enquote{harmonic potential} cases which obscures the covergence orders. However, choosing a high-accuracy integration method is nonetheless very important there as well, as the presence of large gradients (or large directional derivatives) will dominate the exploration.

\subsection{Application: Loss Manifold Analysis for MNIST dataset}\label{sec:application-mnist}

We now consider the MNIST training data set of \num{70000} grey-scale images of 28x28 pixels, see \cite{mnist}. We have divided these into a validation set of \num{10000} images, a test set of \num{5000} images and a training data set of \num{55000} images.

The simplest network to tackle this classification problem is a single layer perceptron with 784 input nodes and 10 output nodes. We use linear output activation and the softmax cross-entropy as loss function. The accuracy is quantified by the strongest output (argmax).

%\subsubsection{Optimization}

%Let us begin with a typical optimization run from the point of view of the data scientist who is given a new unknown dataset\todo{This needs some rewriting and backing by literature about this typical approach.}: Initially, we vary the network parameters such as activation functions or the loss. Eventually, we settle for what seems to be a good parameter set along with suitable batch size and learning rate. Here, we use the softmax cross-entropy as loss function and linear output activations. 
In a first experiment, we use \glostext{SGD} as optimizer with a batch size of 550, i.\,e.~1~\%{} of the dataset size, for 9000 steps and an initial learning rate of 0.5 that is reduced to 0.05 while the batch size is increased to \num{5500} after \num{3000} steps and finally down to 0.01 with no more mini-batching after \num{6000} steps. With this particular training scheme, we obtain a loss of 0.265 and 92.6\% accuracy on the training dataset and a loss of 0.271 and 92.6\% accuracy on the test dataset, c.\,f. 91.6\% and 92.4\% in \cite{LeCun1998}.

\begin{figure}[htbp]
	\centering
	\subfigure[Loss (solid) and accuracy (dotted) along the optimization trajectory consisting of three parts (blue, gray and black) where only every 10th step is shown. The optimization yields 92.6\% test and training set accuracy.]{
		\label{fig:mnist_optimization-step_loss}
		\resizebox{0.48\textwidth}{!}{
			\tikzsetnextfilename{optimization_step-loss_accuracy}
			\ifdefined\myexperiment\empty\else\def\myexperiment{Optimizing_MNIST}\fi%
\ifdefined\mypath\empty\else\def\mypath{\myexperimentspath/13_EnsemblePreconditioning_MNIST/data/\myexperiment}\fi% %

% parameter
\ifdefined\mybeta\empty\else\def\mybeta{10}\fi%
\ifdefined\mydeltat\empty\else\def\mydeltat{0_125}\fi%
\ifdefined\mygamma\empty\else\def\mygamma{10}\fi%
\ifdefined\mysteps\empty\else\def\mysteps{3000}\fi%

\ifdefined\mydof\empty\else\def\mydof{7850}\fi% 784*10+10 (weights + biases)

\pgfplotsset{compat=1.5.1}%

\begin{tikzpicture}[]
\begin{semilogyaxis}[
	width=8cm,height=7cm,
	xmajorgrids,
	scale only axis,
	axis y line*=left,
	xlabel={Step},
	ylabel={Loss},
	legend pos=south east,
	log ticks with fixed point,
	%each nth point=10,
%	restrict y to domain*={50:60},
%	restrict x to domain={0:2e5},
	cycle multi list={
		solid, dashed\nextlist
		blue,gray,black
	},
]
\foreach \mynr/\myleg in {0/first,1/second,2/third}{%
% loss
\addplot+[thick,no marks,black] table[x expr={\mysteps*\mynr+\thisrow{step}},y=loss, col sep=comma]{\mypath/trajectory-train_\myleg-.csv};
}%
\end{semilogyaxis}
\begin{axis}[
	width=8cm,height=7cm,
	xmajorgrids,
	scale only axis,
	axis y line*=right,
	axis x line=none,
	ylabel={Accuracy},
	legend pos=south east,
	%each nth point=10,
%	restrict y to domain*={50:60},
%	restrict x to domain={0:2e5},
	cycle multi list={
		solid, dashed\nextlist
		blue,gray,black
	},
]
\foreach \mynr/\myleg in {0/first,1/second,2/third}{%
% accuracy
\addplot+[thick,dotted, no marks,black] table[x expr={\mysteps*\mynr+\thisrow{step}},y=accuracy, col sep=comma]{\mypath/trajectory-test_\myleg-.csv};
}%
\end{axis}
\end{tikzpicture}
		}
	}
	\subfigure[Loss values along the sampled trajectories where only every 10th step is shown where trajectories obviously deviate rapidly from the initially same starting position. Each run only differs by the random number seed, hence being subject to different thermal noise.]{
		\label{fig:mnist_sampling-step_loss}
		\resizebox{0.46\textwidth}{!}{
			\tikzsetnextfilename{MNIST-sampling_step-loss}
			\input{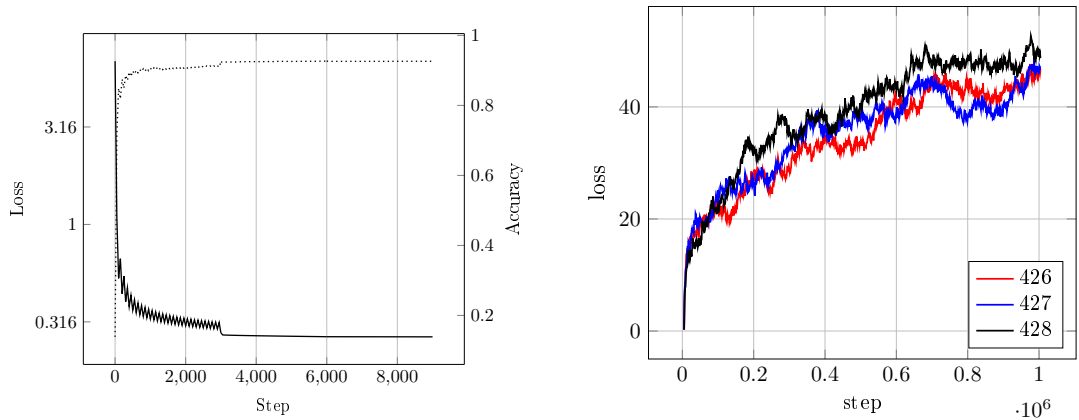}
		}
	}
	\caption{Loss values along the optimization and three sampling trajectories.}
	\label{fig:mnist_trajectories}
\end{figure}

The resulting loss per step is given in Figure~\ref{fig:mnist_optimization-step_loss}. There are fluctuations due to the stochastic gradients that decrease with the learning rate and with the batch size.

%If we stop our analysis here, we might be tempted to believe that we had encountered the global minimum or at least a minimum close to it. At least, the data scientist might be tempted to conclude, it seems good enough for the chosen simple network. 

%\subsubsection{Visualizing the Loss Manifold}

Note that we stop the optimization after a finite number of steps and not when the gradients has a zero norm. Therefore, the points encountered in the loss landscape during optimization are just critical points. Nonetheless, we will refer to them as "quasi-minima" in the following as the gradient norm is very small.

In order to inspect the quality of the quasi-minimum found, we need to look at the resulting loss manifold in a neighborhood. As the network has 7850 degrees of freedom in total, we need appropriate techniques for its visualization. Several of these are discussed in \cite{Li2017}. The current state-of-the-art is to project onto two random directions, as proposed in \cite{Goodfellow2014a}, that frequently only shows little variation~(\cite{Li2017}).

With the sampling approach proposed here, there is a different alternative. Starting in at or near a quasi-minimum, a walker generates a cloud of points iteratively by following the chosen dynamics. Therefore, sampling provides us with a cloud of points that expands within that basin according to rules of these dynamics: if the basin is flat in certain directions, the cloud's expansion will prefer these over directions where the walls are steep. The strength of the preference is controlled by the (inverse) temperature parameter $\beta$ that allows the method to overcome walls to a certain steepness and height.

If we analyze the sampled point cloud's principal components, using the eigensystem of the covariance matrix, then it's major principal component points along the direction where the basin is flattest.

In the following, we compare two components associated with a) the largest eigenvalue, and b) an arbitrarily chosen small eigenvalue of the covarance matrix. This will allow us to get a notion of the extent of the flat part of the basin.

The covariance matrix is computed from a single sampling run of $10^6$ steps using BAOAB as the sampler with a step size $\varepsilon =  0.125$, $\gamma = 10$, $\beta = 10$, and a batch size of 550. 

In the following we have performed three of these sampling runs starting from the same position, equilibrated with \glostext{SGD} with the same batch size of 550, for 5000 steps with a learning rate of 0.5. The only point in using multiple runs is to exclude the possibility that subsequent results depend simply on a specific, rare combination of samples taken. Therefore, each run differs only  by the random number seed and is thus only subject to different thermal noise. Note that the specific initial position will not have an effect as the temperature is high enough to bring the walker quickly to a completely different position: see Figure~\ref{fig:mnist_sampling-step_loss} for loss values of all trajectories.  There, the loss increases from its initial value because of the small value chosen for $\beta$, i.\,e.~a high sampling temperature. The walker typically assumes positions far away from equilibrium.

%We have refrained from employing \glostext{HMC} in this example as it requires exact gradients and therefore has a much higher computational cost.

For reasons of computational efficiency we note the following: to obtain the covariance matrix of very high-dimensional networks the trajectories from several (parallel) runs can be combined to overcome the computational burden.  Moreover, a truncated eigendecomposition, e.\,g.,~using the power method and shift-and-invert, would fully suffice to obtain the largest and a small eigenvalue within a certain range and their associated eigenvectors, taking advantage of symmetry and positive semi-definiteness. 

\begin{figure}[htbp]
	\centering
	%\subfigure[Eigenvalues of the resulting covariance matrices]{
		%\label{fig:mnist_sampling-covariance_eigenvalues}
		\resizebox{0.6\textwidth}{!}{
			\tikzsetnextfilename{MNIST-covariance-eigenvalues}
			\ifdefined\myexperimentname\empty\else\def\myexperimentname{Sampling MNIST}\fi%
\ifdefined\mypath\empty\else\def\mypath{\myexperimentspath/13_EnsemblePreconditioning_MNIST/data/Sampling_MNIST}\fi%

\ifdefined\myseed\empty\else\def\myseed{426}\fi%
\ifdefined\mybatchsize\empty\else\def\mybatchsize{550}\fi%

\def\myMPsigmasq{1.}%
\def\myMPlambda{0.2}%

\begin{tikzpicture}[]
\begin{semilogyaxis}[
	width=8.5cm,height=7.25cm,
	xlabel={Index},
	ylabel={Eigenvalue},
	restrict y to domain={-5:15},
]
\addlegendimage{empty legend}
\addlegendentry{seed}
\foreach \myseed in {426,...,428}{%
\addplot+[thick, no marks] table[x expr={\coordindex}, y=value, col sep=comma]{\mypath/covariance-eigenvalues-batch_size_\mybatchsize-run_seed_\myseed.csv};
\addlegendentryexpanded{\myseed}
}%
\end{semilogyaxis}
\end{tikzpicture}
		}
	%}
	\caption{Logarithmic depiction of the eigenvalues of the covariance matrix of three sampled trajectories. Each legend entry refers to the (different) random number seed employed. Spectra are in good agreement with about 1\% deviation between the three random number seeds.}
	\label{fig:mnist_sampling-covariance_eigenvalues}
\end{figure}
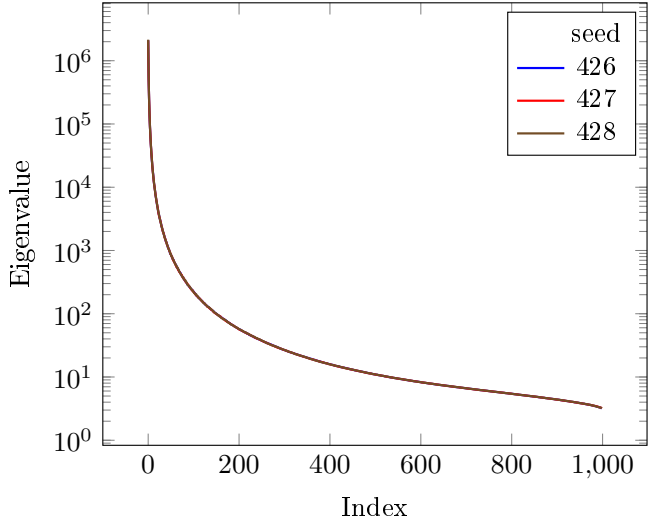

The resulting computed eigenvalue spectra are given in Figure~\ref{fig:mnist_sampling-covariance_eigenvalues}. Note that we have truncated the spectrum here to 1000 non-zero eigenvalues.

All three curves in Figure \ref{fig:mnist_sampling-covariance_eigenvalues} coincide approximately in the logarithmic scale, with a deviation of about 1\%. This leads us to assume that the sampling runs have indeed been long enough for the system to lose memory of its initialization and have produced representative thermodynamic samples from the target ensemble. We hypothesize that this eigenvalue spectrum is representative of the general covariance structure at large scale  for MNIST and for the particular model chosen.

Observe that there is a strong decay of the eigenvalue magnitudes\footnote{The decay resembles the Marchenko-Pastur distribution which describes the singular values of random matrices. %Although, there is some resemblance (in the semi log-scale), we could not obtain a parameter set yielding a good match and 
However, we do not pursue this further in the scope of this article.}.

The decay in the eigenvalue spectrum expresses itself as few eigenvectors with very large eigenvalues and many eigenvectors whose eigenvalues have at least 3 or 4 orders of magnitude smaller eigenvalues. We judge this variation in scale as the spectrum's major feature.

To highlight the extent of the flat part, we pick the first and 64th eigenvalues and use their associated eigenvectors as the directions $v_1$ and $v_0$ in which we plot the loss manifold. The direction $v_1$ represents the large covariance, while $v_0$ represents the small covariance.\footnote{Picking the 1st and the 100th eigenvalue would have given an essentially equivalent visulization.}

The strong decay indicates that random directions are poorly suited for the loss visualization (at least for MNIST and  for the chosen model) as, on average, these will not relate to directions of strong covariance. Hence, the flattest manifold parts will be missed on average when choosing random directions.

For visualizing the loss manifold, we sample it on an equidistant grid with 41 samples per axis along the two chosen directions $v_0 = (W_0, b_0)$, $v_1=(W_1,b_1)$, split into weight components $W_j \in \reel^{784 \times 10}$ and bias components $b_j \in \reel^{10}$. We evaluate $\tilde{L}(c_0,c_1) = \sum_{i=1} \sigma( \tilde{f}(c_0, c_1), y_i)$ with the softmax cross-entropy $\sigma$ and $\tilde{f}(c_0,c_1) = (W + c_0 W_0 + c_1 W_1) \cdot x + B + b_0 + b_1$, i.\,e.~the loss constrained to the two-dimensional subspace and centered at the obtained (local) quasi-minimum at $(W,B)$. We use different intervals per axis of $[-10^i, 10^i]$ with $i \in \{1,0,\ldots,-4,-5\}$, endpoints included.
For visualizing the previously obtained optimization trajectories, we re-evaluate them using the full training dataset (i.\,e. no mini-batches) per step and project them onto the two chosen directions.
Note that also each sampled point of the loss manifold is evaluated using the full training dataset.

%Naturally all of these experiments have been easily possible with the tools available in \emph{TATi}. The covariance matrix itself is readily computed and analysed using functions in the \Pyth{} packages \emph{numpy} and \emph{scipy}.

\begin{figure}[htbp]
	\centering
	\subfigure[${[}-10^{1}, 10^{1}{]}$]{
		\def\mylength{10}%
		\def\mydomain{10}%
		\def\mylevels{10}%
		\def\mysubdomain{\fpeval{\mydomain/10.}}%
		\label{fig:mnist-loss_manifold-i_1}
		\resizebox{0.46\textwidth}{!}{%
			\tikzsetnextfilename{MNIST-loss-manifold-i_1}
			\ifdefined\myexperimentname\empty\else\def\myexperimentname{Sampling minimum basin - covariance}\fi%
\ifdefined\mypath\empty\else\def\mypath{\myexperimentspath/13_EnsemblePreconditioning_MNIST/data/Sampling_minimum_basin_-_covariance}\fi%

\ifdefined\myviewphi\empty\else\def\myviewphi{130}\fi%
\ifdefined\myviewtheta\empty\else\def\myviewtheta{20}\fi%

\ifdefined\mylength\empty\else\def\mylength{0_1}\fi%
\ifdefined\mydomain\empty\else\def\mydomain{0.1}\fi%
\def\minloss{2.64964789e-01}%

\ifdefined\mylevels\empty\else\def\mylevels{10}\fi%

\begin{tikzpicture}[]
% parse last entry as offset
\pgfplotstableread[col sep=comma]{\mypath/../Optimizing_MNIST/trajectory-train_third-.csv}\mytrajectory%

\pgfplotstablegetrowsof\mytrajectory%
\pgfmathparse{\pgfplotsretval-1}%
\xdef\mylastrow{\pgfmathresult}%
\pgfplotstablegetelem{\mylastrow}{c0}\of\mytrajectory%
\xdef\myoptczero{\pgfplotsretval}%
\pgfplotstablegetelem{\mylastrow}{c1}\of\mytrajectory%
\xdef\myoptcone{\pgfplotsretval}%

\begin{axis}[
	width=8.5cm,height=7.25cm,
	xlabel={$c_0$},
	ylabel={$c_1$},
	zlabel={$\ln{|\tilde{L}(c_0,c_1) - L_0|}$},
	%enlargelimits=false,%don't enlarge for filtered z coordinate
	3d box=complete,
	legend style={at={(0.1,1.2)}},
 	view={\myviewphi}{\myviewtheta},
	restrict x to domain={-\mydomain:\mydomain},
	restrict y to domain={-\mydomain:\mydomain},
	zmin=-19,
	%restrict z to domain={0.:10.},
%	colorbar,
	cycle multi list={
		{mark=+,black},{mark=star,blue}, {mark=times,gray}
	},
	point meta rel=per plot,
]

% CONTOUR
\addplot3[%
	solid,thick,
	mesh/rows=41,mesh/cols=41,mesh/ordering=y varies,
	contour gnuplot={
		contour dir=z,
		labels=false,
		draw color=black,
		number={\mylevels},
		%levels={1e-7,1e-6,1e-5,1e-4},
	},
	%point meta min=0,
	%point meta max=1,
	z filter/.code={\def\pgfmathresult{-19}},%
]%
	table[
		x=c0,y=c1,z expr={((\thisrow{loss}-\minloss))},
		col sep=comma]
{\mypath/samples_covariance-interval_length_\mylength.csv};
%\addlegendentry{contour}

% manifold
\addplot3+[%
	no marks,
	mesh,mesh/rows=41,mesh/cols=41,mesh/ordering=y varies,
	surf,
	shader=interp,
	colorbar source,
] table[x=c0,y=c1,z expr={ln((\thisrow{loss}-\minloss))},col sep=comma]{\mypath/samples_covariance-interval_length_\mylength.csv};
%\addlegendentry{loss manifold}

% trajectory
\pgfplotsinvokeforeach{first,second,third}{%
\addplot3+[%
	%only marks,
	%mark=+,
	each nth point=1,
	restrict expr to domain={\thisrow{loss}}{2.64964789e-01:1},
	unbounded coords=discard,
] table[x expr={\thisrow{c0}-\myoptczero},y expr={\thisrow{c1}-\myoptcone},z expr={ln((\thisrow{loss}-\minloss))}, col sep=comma]{\mypath/../Optimizing_MNIST/trajectory-train_#1-.csv};
%\addlegendentry{#1 opt run}
}
\ifdefined\mysubdomain{
\draw[fill=white,draw=red,thick] (axis cs:-\mysubdomain,-\mysubdomain,-19) -- (axis cs:\mysubdomain,-\mysubdomain,-19) -- (axis cs:\mysubdomain,\mysubdomain,-19) -- (axis cs:-\mysubdomain,\mysubdomain,-19) -- cycle;}%
\fi
\end{axis}
\end{tikzpicture}
		}
	}
	\subfigure[${[}-10^{0}, 10^{0}{]}$]{
		\def\mylength{1}%
		\def\mydomain{1}%
		\def\mylevels{10}%
		\def\mysubdomain{\fpeval{\mydomain/10.}}%
		\label{fig:mnist-loss_manifold-i_0}
		\resizebox{0.46\textwidth}{!}{%
			\tikzsetnextfilename{MNIST-loss-manifold-i_0}
			\ifdefined\myexperimentname\empty\else\def\myexperimentname{Sampling minimum basin - covariance}\fi%
\ifdefined\mypath\empty\else\def\mypath{\myexperimentspath/13_EnsemblePreconditioning_MNIST/data/Sampling_minimum_basin_-_covariance}\fi%

\ifdefined\myviewphi\empty\else\def\myviewphi{130}\fi%
\ifdefined\myviewtheta\empty\else\def\myviewtheta{20}\fi%

\ifdefined\mylength\empty\else\def\mylength{0_1}\fi%
\ifdefined\mydomain\empty\else\def\mydomain{0.1}\fi%
\def\minloss{2.64964789e-01}%

\ifdefined\mylevels\empty\else\def\mylevels{10}\fi%

\begin{tikzpicture}[]
% parse last entry as offset
\pgfplotstableread[col sep=comma]{\mypath/../Optimizing_MNIST/trajectory-train_third-.csv}\mytrajectory%

\pgfplotstablegetrowsof\mytrajectory%
\pgfmathparse{\pgfplotsretval-1}%
\xdef\mylastrow{\pgfmathresult}%
\pgfplotstablegetelem{\mylastrow}{c0}\of\mytrajectory%
\xdef\myoptczero{\pgfplotsretval}%
\pgfplotstablegetelem{\mylastrow}{c1}\of\mytrajectory%
\xdef\myoptcone{\pgfplotsretval}%

\begin{axis}[
	width=8.5cm,height=7.25cm,
	xlabel={$c_0$},
	ylabel={$c_1$},
	zlabel={$\ln{|\tilde{L}(c_0,c_1) - L_0|}$},
	%enlargelimits=false,%don't enlarge for filtered z coordinate
	3d box=complete,
	legend style={at={(0.1,1.2)}},
 	view={\myviewphi}{\myviewtheta},
	restrict x to domain={-\mydomain:\mydomain},
	restrict y to domain={-\mydomain:\mydomain},
	zmin=-19,
	%restrict z to domain={0.:10.},
%	colorbar,
	cycle multi list={
		{mark=+,black},{mark=star,blue}, {mark=times,gray}
	},
	point meta rel=per plot,
]

% CONTOUR
\addplot3[%
	solid,thick,
	mesh/rows=41,mesh/cols=41,mesh/ordering=y varies,
	contour gnuplot={
		contour dir=z,
		labels=false,
		draw color=black,
		number={\mylevels},
		%levels={1e-7,1e-6,1e-5,1e-4},
	},
	%point meta min=0,
	%point meta max=1,
	z filter/.code={\def\pgfmathresult{-19}},%
]%
	table[
		x=c0,y=c1,z expr={((\thisrow{loss}-\minloss))},
		col sep=comma]
{\mypath/samples_covariance-interval_length_\mylength.csv};
%\addlegendentry{contour}

% manifold
\addplot3+[%
	no marks,
	mesh,mesh/rows=41,mesh/cols=41,mesh/ordering=y varies,
	surf,
	shader=interp,
	colorbar source,
] table[x=c0,y=c1,z expr={ln((\thisrow{loss}-\minloss))},col sep=comma]{\mypath/samples_covariance-interval_length_\mylength.csv};
%\addlegendentry{loss manifold}

% trajectory
\pgfplotsinvokeforeach{first,second,third}{%
\addplot3+[%
	%only marks,
	%mark=+,
	each nth point=1,
	restrict expr to domain={\thisrow{loss}}{2.64964789e-01:1},
	unbounded coords=discard,
] table[x expr={\thisrow{c0}-\myoptczero},y expr={\thisrow{c1}-\myoptcone},z expr={ln((\thisrow{loss}-\minloss))}, col sep=comma]{\mypath/../Optimizing_MNIST/trajectory-train_#1-.csv};
%\addlegendentry{#1 opt run}
}
\ifdefined\mysubdomain{
\draw[fill=white,draw=red,thick] (axis cs:-\mysubdomain,-\mysubdomain,-19) -- (axis cs:\mysubdomain,-\mysubdomain,-19) -- (axis cs:\mysubdomain,\mysubdomain,-19) -- (axis cs:-\mysubdomain,\mysubdomain,-19) -- cycle;}%
\fi
\end{axis}
\end{tikzpicture}
		}
	}
	\subfigure[${[}-10^{-1}, 10^{-1}{]}$]{
		\def\mylength{0_1}%
		\def\mydomain{0.1}%
		\def\mylevels{20}%
		\def\mysubdomain{\fpeval{\mydomain/10.}}%
		\label{fig:mnist-loss_manifold-i_-1}
		\resizebox{0.46\textwidth}{!}{%
			\tikzsetnextfilename{MNIST-loss-manifold-i_-1}
			\ifdefined\myexperimentname\empty\else\def\myexperimentname{Sampling minimum basin - covariance}\fi%
\ifdefined\mypath\empty\else\def\mypath{\myexperimentspath/13_EnsemblePreconditioning_MNIST/data/Sampling_minimum_basin_-_covariance}\fi%

\ifdefined\myviewphi\empty\else\def\myviewphi{130}\fi%
\ifdefined\myviewtheta\empty\else\def\myviewtheta{20}\fi%

\ifdefined\mylength\empty\else\def\mylength{0_1}\fi%
\ifdefined\mydomain\empty\else\def\mydomain{0.1}\fi%
\def\minloss{2.64964789e-01}%

\ifdefined\mylevels\empty\else\def\mylevels{10}\fi%

\begin{tikzpicture}[]
% parse last entry as offset
\pgfplotstableread[col sep=comma]{\mypath/../Optimizing_MNIST/trajectory-train_third-.csv}\mytrajectory%

\pgfplotstablegetrowsof\mytrajectory%
\pgfmathparse{\pgfplotsretval-1}%
\xdef\mylastrow{\pgfmathresult}%
\pgfplotstablegetelem{\mylastrow}{c0}\of\mytrajectory%
\xdef\myoptczero{\pgfplotsretval}%
\pgfplotstablegetelem{\mylastrow}{c1}\of\mytrajectory%
\xdef\myoptcone{\pgfplotsretval}%

\begin{axis}[
	width=8.5cm,height=7.25cm,
	xlabel={$c_0$},
	ylabel={$c_1$},
	zlabel={$\ln{|\tilde{L}(c_0,c_1) - L_0|}$},
	%enlargelimits=false,%don't enlarge for filtered z coordinate
	3d box=complete,
	legend style={at={(0.1,1.2)}},
 	view={\myviewphi}{\myviewtheta},
	restrict x to domain={-\mydomain:\mydomain},
	restrict y to domain={-\mydomain:\mydomain},
	zmin=-19,
	%restrict z to domain={0.:10.},
%	colorbar,
	cycle multi list={
		{mark=+,black},{mark=star,blue}, {mark=times,gray}
	},
	point meta rel=per plot,
]

% CONTOUR
\addplot3[%
	solid,thick,
	mesh/rows=41,mesh/cols=41,mesh/ordering=y varies,
	contour gnuplot={
		contour dir=z,
		labels=false,
		draw color=black,
		number={\mylevels},
		%levels={1e-7,1e-6,1e-5,1e-4},
	},
	%point meta min=0,
	%point meta max=1,
	z filter/.code={\def\pgfmathresult{-19}},%
]%
	table[
		x=c0,y=c1,z expr={((\thisrow{loss}-\minloss))},
		col sep=comma]
{\mypath/samples_covariance-interval_length_\mylength.csv};
%\addlegendentry{contour}

% manifold
\addplot3+[%
	no marks,
	mesh,mesh/rows=41,mesh/cols=41,mesh/ordering=y varies,
	surf,
	shader=interp,
	colorbar source,
] table[x=c0,y=c1,z expr={ln((\thisrow{loss}-\minloss))},col sep=comma]{\mypath/samples_covariance-interval_length_\mylength.csv};
%\addlegendentry{loss manifold}

% trajectory
\pgfplotsinvokeforeach{first,second,third}{%
\addplot3+[%
	%only marks,
	%mark=+,
	each nth point=1,
	restrict expr to domain={\thisrow{loss}}{2.64964789e-01:1},
	unbounded coords=discard,
] table[x expr={\thisrow{c0}-\myoptczero},y expr={\thisrow{c1}-\myoptcone},z expr={ln((\thisrow{loss}-\minloss))}, col sep=comma]{\mypath/../Optimizing_MNIST/trajectory-train_#1-.csv};
%\addlegendentry{#1 opt run}
}
\ifdefined\mysubdomain{
\draw[fill=white,draw=red,thick] (axis cs:-\mysubdomain,-\mysubdomain,-19) -- (axis cs:\mysubdomain,-\mysubdomain,-19) -- (axis cs:\mysubdomain,\mysubdomain,-19) -- (axis cs:-\mysubdomain,\mysubdomain,-19) -- cycle;}%
\fi
\end{axis}
\end{tikzpicture}
		}
	}
	\subfigure[${[}-10^{-2}, 10^{-2}{]}$]{
		\def\mylength{0_01}%
		\def\mydomain{0.01}%
		\def\mylevels{4}%
		\label{fig:mnist-loss_manifold-i_-2}
		\resizebox{0.46\textwidth}{!}{%
			\tikzsetnextfilename{MNIST-loss-manifold-i_-2}
			\ifdefined\myexperimentname\empty\else\def\myexperimentname{Sampling minimum basin - covariance}\fi%
\ifdefined\mypath\empty\else\def\mypath{\myexperimentspath/13_EnsemblePreconditioning_MNIST/data/Sampling_minimum_basin_-_covariance}\fi%

\ifdefined\myviewphi\empty\else\def\myviewphi{130}\fi%
\ifdefined\myviewtheta\empty\else\def\myviewtheta{20}\fi%

\ifdefined\mylength\empty\else\def\mylength{0_1}\fi%
\ifdefined\mydomain\empty\else\def\mydomain{0.1}\fi%
\def\minloss{2.64964789e-01}%

\ifdefined\mylevels\empty\else\def\mylevels{10}\fi%

\begin{tikzpicture}[]
% parse last entry as offset
\pgfplotstableread[col sep=comma]{\mypath/../Optimizing_MNIST/trajectory-train_third-.csv}\mytrajectory%

\pgfplotstablegetrowsof\mytrajectory%
\pgfmathparse{\pgfplotsretval-1}%
\xdef\mylastrow{\pgfmathresult}%
\pgfplotstablegetelem{\mylastrow}{c0}\of\mytrajectory%
\xdef\myoptczero{\pgfplotsretval}%
\pgfplotstablegetelem{\mylastrow}{c1}\of\mytrajectory%
\xdef\myoptcone{\pgfplotsretval}%

\begin{axis}[
	width=8.5cm,height=7.25cm,
	xlabel={$c_0$},
	ylabel={$c_1$},
	zlabel={$\ln{|\tilde{L}(c_0,c_1) - L_0|}$},
	%enlargelimits=false,%don't enlarge for filtered z coordinate
	3d box=complete,
	legend style={at={(0.1,1.2)}},
 	view={\myviewphi}{\myviewtheta},
	restrict x to domain={-\mydomain:\mydomain},
	restrict y to domain={-\mydomain:\mydomain},
	zmin=-19,
	%restrict z to domain={0.:10.},
%	colorbar,
	cycle multi list={
		{mark=+,black},{mark=star,blue}, {mark=times,gray}
	},
	point meta rel=per plot,
]

% CONTOUR
\addplot3[%
	solid,thick,
	mesh/rows=41,mesh/cols=41,mesh/ordering=y varies,
	contour gnuplot={
		contour dir=z,
		labels=false,
		draw color=black,
		number={\mylevels},
		%levels={1e-7,1e-6,1e-5,1e-4},
	},
	%point meta min=0,
	%point meta max=1,
	z filter/.code={\def\pgfmathresult{-19}},%
]%
	table[
		x=c0,y=c1,z expr={((\thisrow{loss}-\minloss))},
		col sep=comma]
{\mypath/samples_covariance-interval_length_\mylength.csv};
%\addlegendentry{contour}

% manifold
\addplot3+[%
	no marks,
	mesh,mesh/rows=41,mesh/cols=41,mesh/ordering=y varies,
	surf,
	shader=interp,
	colorbar source,
] table[x=c0,y=c1,z expr={ln((\thisrow{loss}-\minloss))},col sep=comma]{\mypath/samples_covariance-interval_length_\mylength.csv};
%\addlegendentry{loss manifold}

% trajectory
\pgfplotsinvokeforeach{first,second,third}{%
\addplot3+[%
	%only marks,
	%mark=+,
	each nth point=1,
	restrict expr to domain={\thisrow{loss}}{2.64964789e-01:1},
	unbounded coords=discard,
] table[x expr={\thisrow{c0}-\myoptczero},y expr={\thisrow{c1}-\myoptcone},z expr={ln((\thisrow{loss}-\minloss))}, col sep=comma]{\mypath/../Optimizing_MNIST/trajectory-train_#1-.csv};
%\addlegendentry{#1 opt run}
}
\ifdefined\mysubdomain{
\draw[fill=white,draw=red,thick] (axis cs:-\mysubdomain,-\mysubdomain,-19) -- (axis cs:\mysubdomain,-\mysubdomain,-19) -- (axis cs:\mysubdomain,\mysubdomain,-19) -- (axis cs:-\mysubdomain,\mysubdomain,-19) -- cycle;}%
\fi
\end{axis}
\end{tikzpicture}
		}
	}
	\caption{Visualization of the regularly sampled loss manifold for the MNIST dataset of single-layer perceptron with softmax cross-entropy loss function and linear activations. The two coordinate directions $c_0$, $c_1$ correspond to the 1st and 64th eigenvalues (descending) of a covariance matrix sampled from a very long run. The $z$-axis, in \emph{log scale}, corresponds to $\ln{(\tilde{L}(c_0,c_1)-L_0)}$, where $L_0$ is the smallest loss encountered overall. The value for $z$ is also used to color the manifold. Additionally, an optimization trajectory using \protect\glostext{SGD} consisting of three consecutive parts is given where after each leg the learning rate is reduced and batch size is increased. The terminal point of its last leg provides the point of origin. Only every 50th step is shown. Note that we show the exact loss for the whole training dataset for both manifold and trajectory. A red square shows the subsequent plot's domain.}
	\label{fig:mnist-loss_manifold}
\end{figure}
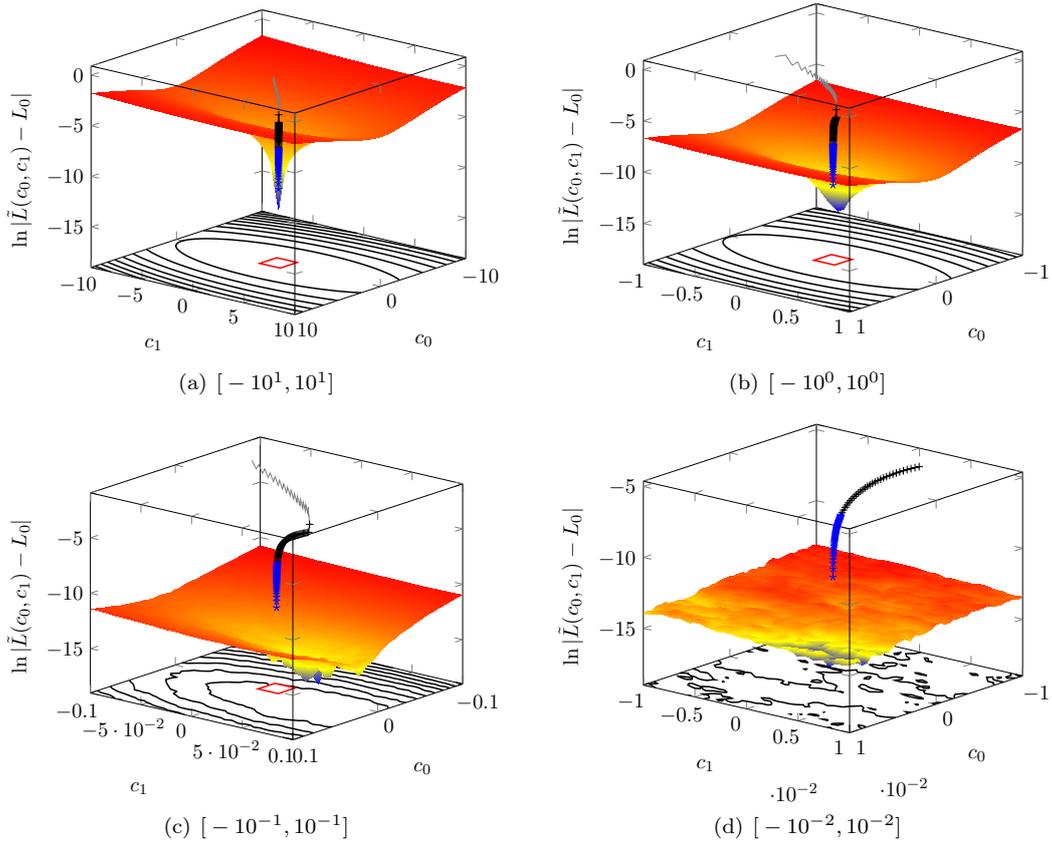

In Figure~\ref{fig:mnist-loss_manifold} we look at four of these manifold plots where we give the sampled manifold and the projected optimization trajectory. The $x$ and $y$ axes correspond to the two chosen directions, the $z$ axis gives the loss as $\ln{ \mid \tilde{L}(c_0,c_1)-L_0) \mid }$, with respect to the lowest loss $L_0$ value found overall, at the specific point $(c_0,c_10)$ in case of the sampled manifold and the true loss in case of the projected trajectory.  Because of the projection the trajectory steps will not lie on the manifold itself.

In Figure~\ref{fig:mnist-loss_manifold-i_1} we recognize a large funnel where the $c_0$ direction is associated with the large eigenvalue and the $c_1$ direction is associated with the small eigenvalue. We see that the optimization trajectory gradually enters the funnel in \subref{fig:mnist-loss_manifold-i_0}{}. However, in \subref{fig:mnist-loss_manifold-i_-1}{} we realize that the previously obtained optimization trajectory ends prematurely, not in the possibly global minimum but stuck in one of the lower minima on the funnel wall in \subref{fig:mnist-loss_manifold-i_-2}{}.

%We note without showing further visualizations that this general structure is valid also for other eigenvector pairs picked from the covariance matrix. The extent of the anisotropy of the funnel depends on the ratio of the associated eigenvalues of the two chosen directions.\todo{Re-check that this is indeed true.}

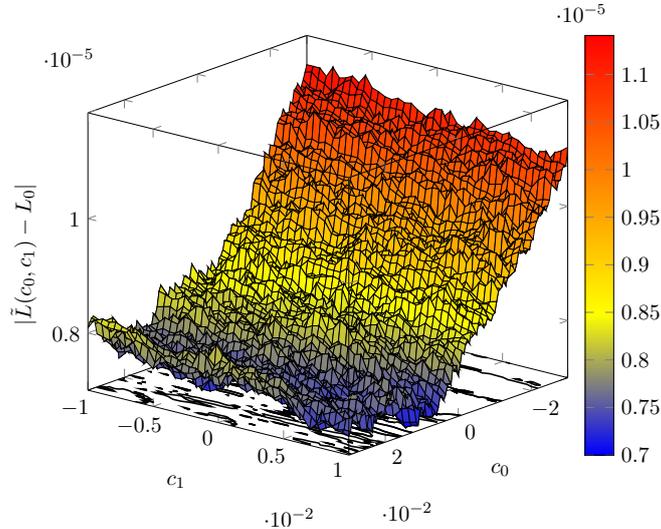
\begin{figure}[htbp]
	\centering
	\resizebox{0.6\textwidth}{!}{%
		\def\mylength{0_01}%
		\def\mydomain{0.01}%
		\def\mylevels{3}%
		\tikzsetnextfilename{MNIST-loss-manifold-normal_scale-stripe}
		\ifdefined\myexperimentname\empty\else\def\myexperimentname{Sampling minimum basin - covariance}\fi%
\ifdefined\mypath\empty\else\def\mypath{\myexperimentspath/17_MNIST_Loss_Manifold/data/Scanning_Loss_Manifold}\fi%

\ifdefined\myviewphi\empty\else\def\myviewphi{130}\fi%
\ifdefined\myviewtheta\empty\else\def\myviewtheta{20}\fi%

\ifdefined\mylength\empty\else\def\mylength{0_01}\fi%
\ifdefined\mydomain\empty\else\def\mydomain{0.01}\fi%
\def\minloss{2.64964789e-01}%

\ifdefined\myoffsetstep\empty\else\def\myoffsetstep{0}\fi%
\ifdefined\myoffsety\empty\else\def\myoffsety{0_}\fi%

\ifdefined\mylevels\empty\else\def\mylevels{3}\fi%

\xdef\mysubdomain{\fpeval{\mydomain/10.}}%
\begin{tikzpicture}[]

% parse last entry as offset
\pgfplotstableread[col sep=comma]{\mypath/../../../13_EnsemblePreconditioning_MNIST/data/Optimizing_MNIST/trajectory-train_third-.csv}\mytrajectory%

\pgfplotstablegetrowsof\mytrajectory%
\pgfmathparse{\pgfplotsretval-1}%
\xdef\mylastrow{\pgfmathresult}%
\pgfplotstablegetelem{\mylastrow}{c0}\of\mytrajectory%
\xdef\myoptczero{\pgfplotsretval}%
\pgfplotstablegetelem{\mylastrow}{c1}\of\mytrajectory%
\xdef\myoptcone{\pgfplotsretval}%

\begin{axis}[
	width=10cm,height=9cm,
	xlabel={$c_0$},
	ylabel={$c_1$},
	zlabel={$|\tilde{L}(c_0,c_1) - L_0|$},
	%enlargelimits=false,%don't enlarge for filtered z coordinate
	3d box=complete,
	legend style={at={(0.1,1.2)}},
 	view={\myviewphi}{\myviewtheta},
	%restrict x to domain={-\mydomain:\mydomain},
	%restrict y to domain={-\mydomain:\mydomain},
	zmin=7e-6,
	%restrict z to domain={0.:10.},
	cycle multi list={
		black, blue, gray
	},
	%point meta rel=per plot,
	colorbar,
]
% CONTOUR
\foreach \myoffsetstep in {-2,-1,...,2}{%
\addplot3[
	solid,thick,
	mesh/rows=41,mesh/cols=41,mesh/ordering=y varies,
	contour gnuplot={
		contour dir=z,
		labels=false,
		draw color=black,
		number={\mylevels},
		%levels={1e-7,1e-6,1e-5,1e-4},
	},
	%point meta min=0,
	%point meta max=1,
	z filter/.code={\def\pgfmathresult{7e-6}},%
]%
	table[
		x=c0,y=c1,z expr={(abs(\thisrow{loss}-\minloss))},
		col sep=comma]
{\mypath/samples_covariance-interval_length_\mylength-offset_step_\myoffsetstep-offset_y_\myoffsety.csv};
}
%\addlegendentry{contour}

% manifold
\foreach \myoffsetstep in {-2,-1,...,2}{%
\addplot3+[
	no marks,
	mesh,mesh/rows=41,mesh/cols=41,mesh/ordering=y varies,
	surf,
	faceted color=black,
	%shader=interp,
	%colorbar source,
] table[x=c0,y=c1,z expr={(abs(\thisrow{loss}-\minloss))},col sep=comma]{\mypath/samples_covariance-interval_length_\mylength-offset_step_\myoffsetstep-offset_y_\myoffsety.csv};
}%
%\addlegendentry{loss manifold}

% trajectory
% \pgfplotsinvokeforeach{first,second,third}{%
% \addplot3+[
% 	%only marks,
% 	mark=+,
% 	each nth point=1,
% 	restrict expr to domain={\thisrow{loss}}{2.64964789e-01:1},
% 	unbounded coords=discard,
% ] table[x expr={\thisrow{c0}-\myoptczero},y expr={\thisrow{c1}-\myoptcone},z expr={ln((\thisrow{loss}-\minloss))}, col sep=comma]{\mypath/../../../13_EnsemblePreconditioning_MNIST/data/Optimizing_MNIST/trajectory-train_#1-.csv};
% \addlegendentry{#1 opt run}
% }%

%\draw[fill=white,draw=red,thick] (axis cs:-\mysubdomain,-\mysubdomain,-19) -- (axis cs:\mysubdomain,-\mysubdomain,-19) -- (axis cs:\mysubdomain,\mysubdomain,-19) -- (axis cs:-\mysubdomain,\mysubdomain,-19) -- cycle;
\end{axis}
\end{tikzpicture}
	}
	\caption{Sweep in the direction $c_0$ over five times the domain length with respect to Figure~\ref{fig:mnist-loss_manifold-i_-2}. The walls of the funnel's bottom are dented with many local minima. Note that the figure is in linear scale for comparison to the log-scale before.}
	\label{fig:mnist-loss_manifold-sweep}
\end{figure}

The intuition therefore is that the loss manifold resembles a large funnel whose extension though differs significantly in each direction. Its walls are corrugated with many local minima, see Figure~\ref{fig:mnist-loss_manifold-sweep}, especially at its bottom. 
Note that this funnel is not a product of the logarithmic scale, again see Figure~\ref{fig:mnist-loss_manifold-sweep} where a normal scale is used.
This corresponds well with the observation in \cite{Goodfellow2014a} that optimization runs never encounter serious obstacles and that there is a \enquote{sea of minima} in a small band of the loss bounded from below as proposed in \cite{Choromanska2015} from translating results on spherical spin-glasses. Note though that the latter was derived for a different loss function.

This study hints at the usefulness of second-order methods also for optimization, 
see \cite{Nocedal2018} for a review.

Even in the limited scope of the two-dimensional projection we have seen that a minimum associated with an even smaller loss would have been attainable. However, the differences in the loss values are marginal, namely less than $10^{-5}$.
%Note that these minima all have around $10^{-4}$ or less difference to one another.
Naturally, this analysis is not complete. There are possibly multiple funnels from multiple minima with high barriers in between, as in the disconnectivity graphs by \cite{Ballard2017}. Recent work~(\cite{Draxler2018}), however, suggests the contrary and corroborates our finding of a single funnel with many minima at its bottom. We conclude by remarking that this type of analysis can also easily be extended to multi-layer perceptrons and potentially to more advanced network architectures. Note that for multi-layer perceptron, scale invariance needs to be accounted for, see \cite{Li2017} for a normalization scheme. There, multiple minima may be encountered due to symmetries.

\subsection{Application: More complex networks for MNIST dataset}\label{sec:application-mnist-hidden}

\pgfkeys{/pgf/fpu=true}

%\pgfmathparse{round(784*10+10*10+10+10)}%
%\pgfmathprintnumberto{\pgfmathresult}{\numberdofsmall}%

\pgfmathparse{round(784*100+100*10+100+10)}%
\pgfmathprintnumberto{\pgfmathresult}{\numberdoflarge}%
\pgfkeys{/pgf/fpu=false}

The previous model with no hidden layer used linear activations; only its loss function was non-linear. In this subsection, we consider very briefly the same dataset with a multi-layer perceptron with a single hidden layer of 100 nodes and sigmoidal hidden activation functions. This is to illustrate that the technique is not limited to linear models. We use again linear activations in the output layer and softmax cross-entropy as loss functions. This model has \numberdoflarge{} degrees of freedom, i.\,e.~an order of magnitude more than the linear model used before. We will not distinguish in the following between degrees of freedom associated to the first hidden layer and those associated to the output layer.

We have implemented \glostext{GD} with a Barzilei-Borwein step width choice (\glsentrytext{glos:BBGD}), see \cite{Tan2016}, using the full gradient information. We use 3000 \glsentrytext{glos:BBGD} steps with an initial learning rate of 0.05. This advanced optimization method is not required for the sampling's starting point, there \glostext{SGD} would suffice. However, it yields an improved origin of the loss landscape in the following visualization.

Covariance directions have been obtained through sampling with the same parameters as with the linear model. Then, centered around the located quasi-minimum and using the 1st and 100th covariance direction sorted by their associated eigenvalue, we again sample on an equidistant grid in the subspace spanned by these two directions $v_1$ and $v_0$. Again, the sampling still finds parametrizations $\theta$ with smaller loss than located during optimization whose lowest serves as reference $L_0$.

\begin{figure}[htbp]
	% minima extracted using 22_LandscapeAnalysis_MNIST/scripts/extractMinimumLoss.sh.
	\centering
	\def\mydimension{100}%
	\def\mylevels{20}%
	\def\myactivation{sigmoid}%
	\def\minloss{0.000007585034380}%
	\def\myviewphi{130}%
	\subfigure[${[}-10^{1}, 10^{1}{]}$]{
		\def\mylength{10}%
		\def\mydomain{10}%
		\def\myviewtheta{10}%
		\def\myzmin{-15}%
		\label{fig:mnist_complex-loss_manifold-sigmoid-absolute}
		\resizebox{0.42\textwidth}{!}{%
			\tikzsetnextfilename{MNIST_complex-loss-manifold-sigmoid-absolute}
			\ifdefined\myexperimentname\empty\else\def\myexperimentname{Sampling minimum basin - covariance}\fi%
\ifdefined\mypath\empty\else\def\mypath{\myexperimentspath/22_LandscapeAnalysis_MNIST/data/Sampling_minimum_basin_-_covariance}\fi%

\ifdefined\myviewphi\empty\else\def\myviewphi{130}\fi%
\ifdefined\myviewtheta\empty\else\def\myviewtheta{10}\fi%

% parameters
\ifdefined\myactivation\empty\else\def\myactivation{sigmoid}\fi%
\ifdefined\mydimension\empty\else\def\mydimension{100}\fi%
\ifdefined\myloss\empty\else\def\myloss{softmax_cross_entropy}\fi%

\ifdefined\mylength\empty\else\def\mylength{10}\fi%
\ifdefined\mydomain\empty\else\def\mydomain{10}\fi%
\def\minloss{0.}%000007585034380}%

\ifdefined\mylevels\empty\else\def\mylevels{10}\fi%

\xdef\mysubdomain{\fpeval{\mydomain/10.}}%
\ifdefined\myzmin\empty\else\def\myzmin{-15}\fi%

\begin{tikzpicture}[]
\begin{axis}[
	xlabel={c0},
	ylabel={c1},
	zlabel={$ln(\tilde{L}(c_0,c_1))$},
	%enlargelimits=false,%don't enlarge for filtered z coordinate
	3d box=complete,
	legend style={at={(0.1,1.2)}},
 	view={\myviewphi}{\myviewtheta},
	restrict x to domain={-\mydomain:\mydomain},
	restrict y to domain={-\mydomain:\mydomain},
	zmin=\myzmin,
	%restrict z to domain={0.:10.},
	cycle multi list={
		black, blue, gray
	},
	point meta rel=per plot,
	%colorbar,
]
% CONTOUR
\addplot3[
	solid,thick,
	mesh/rows=41,mesh/cols=41,mesh/ordering=y varies,
	contour gnuplot={
		contour dir=z,
		labels=false,
		draw color=black,
		number={\mylevels},
		%levels={1e-7,1e-6,1e-5,1e-4},
	},
	%point meta min=0,
	%point meta max=1,
	z filter/.code={\def\pgfmathresult{\myzmin}},%
]%
	table[
		x=c0,y=c1,z expr={((\thisrow{loss}-\minloss))},
		col sep=comma]
{\mypath/samples_covariance-hidden_activation_\myactivation-hidden_dimension_\mydimension-interval_length_\mylength-loss_\myloss.csv};
%\addlegendentry{contour}

% manifold
\addplot3+[
	no marks,
	mesh,mesh/rows=41,mesh/cols=41,mesh/ordering=y varies,
	surf,
	shader=interp,
	colorbar source,
] table[x=c0,y=c1,z expr={ln((\thisrow{loss}-\minloss))},col sep=comma]{\mypath/samples_covariance-hidden_activation_\myactivation-hidden_dimension_\mydimension-interval_length_\mylength-loss_\myloss.csv};
%\addlegendentry{loss manifold}

%\draw[fill=white,draw=red,thick] (axis cs:-\mysubdomain,-\mysubdomain,\myzmin) -- (axis cs:\mysubdomain,-\mysubdomain,\myzmin) -- (axis cs:\mysubdomain,\mysubdomain,\myzmin) -- (axis cs:-\mysubdomain,\mysubdomain,\myzmin) -- cycle;
\end{axis}
\end{tikzpicture}
		}
	}
	\subfigure[${[}-10^{-1}, 10^{-1}{]}$]{
		\def\mylength{0_1}%
		\def\mydomain{0.1}%
		\def\myviewtheta{20}%
		\def\myzmin{-24}%
		\label{fig:mnist_complex-loss_manifold-sigmoid-global}
		\resizebox{0.42\textwidth}{!}{%
			\tikzsetnextfilename{MNIST_complex-loss-manifold-sigmoid-global}
			\ifdefined\myexperimentname\empty\else\def\myexperimentname{Sampling minimum basin - covariance}\fi%
\ifdefined\mypath\empty\else\def\mypath{\myexperimentspath/22_LandscapeAnalysis_MNIST/data/Sampling_minimum_basin_-_covariance}\fi%

\ifdefined\myviewphi\empty\else\def\myviewphi{130}\fi%
\ifdefined\myviewtheta\empty\else\def\myviewtheta{80}\fi%

% parameters
\ifdefined\myactivation\empty\else\def\myactivation{sigmoid}\fi%
\ifdefined\mydimension\empty\else\def\mydimension{100}\fi%
\ifdefined\myloss\empty\else\def\myloss{softmax_cross_entropy}\fi%

\ifdefined\mylength\empty\else\def\mylength{0_01}\fi%
\ifdefined\mydomain\empty\else\def\mydomain{0.01}\fi%
\ifdefined\minloss\empty\else\def\minloss{0.000007585034380}\fi%

\ifdefined\mylevels\empty\else\def\mylevels{20}\fi%

\xdef\mysubdomain{\fpeval{\mydomain/10.}}%
\ifdefined\myzmin\empty\else\def\myzmin{-80}\fi%

\begin{tikzpicture}[]
\begin{axis}[
	xlabel={c0},
	ylabel={c1},
	zlabel={$ln(\tilde{L}(c_0,c_1)-L_0)$},
	%enlargelimits=false,%don't enlarge for filtered z coordinate
	3d box=complete,
	legend style={at={(0.1,1.2)}},
 	view={\myviewphi}{\myviewtheta},
	restrict x to domain={-\mydomain:\mydomain},
	restrict y to domain={-\mydomain:\mydomain},
	zmin=\myzmin,
	%restrict z to domain={0.:10.},
	cycle multi list={
		black, blue, gray
	},
	point meta rel=per plot,
	%colorbar,
]
% CONTOUR
\addplot3[
	solid,thick,
	mesh/rows=41,mesh/cols=41,mesh/ordering=y varies,
	contour gnuplot={
		contour dir=z,
		labels=false,
		draw color=black,
		number={\mylevels},
		%levels={1e-7,1e-6,1e-5,1e-4},
	},
	%point meta min=0,
	%point meta max=1,
	z filter/.code={\def\pgfmathresult{\myzmin}},%
]%
	table[
		x=c0,y=c1,z expr={((\thisrow{loss}-\minloss))},
		col sep=comma]
{\mypath/samples_covariance-hidden_activation_\myactivation-hidden_dimension_\mydimension-interval_length_\mylength-loss_\myloss.csv};
%\addlegendentry{contour}

% manifold
\addplot3+[
	no marks,
	mesh,mesh/rows=41,mesh/cols=41,mesh/ordering=y varies,
	surf,
	shader=interp,
	colorbar source,
] table[x=c0,y=c1,z expr={ln((\thisrow{loss}-\minloss))},col sep=comma]{\mypath/samples_covariance-hidden_activation_\myactivation-hidden_dimension_\mydimension-interval_length_\mylength-loss_\myloss.csv};
%\addlegendentry{loss manifold}

\end{axis}
\end{tikzpicture}
		}
	}
	\caption{Visualization of the regularly sampled loss manifold for the MNIST dataset of multi-layer perceptron with 100 hidden nodes, softmax cross-entropy loss function, sigmoid hidden activations and linear output activations. The two coordinate directions $c_0$, $c_1$ correspond to directions associated with the 1st and 100th eigenvalues (descending) of a covariance matrix sampled from a very long run. The value for $z$, also used to color the manifold, is as follows. \subref{fig:mnist_complex-loss_manifold-sigmoid-absolute}: The $z$-axis, in \emph{log scale}, corresponds to $\ln{(\tilde{L}(c_0,c_1))}$ with $c_0,c_1 \in [-10,10]$, i.\,e.~the absolute logarithmic loss on the largest sampled interval in the two chosen covariance directions. \subref{fig:mnist_complex-loss_manifold-sigmoid-global}: The $z$-axis, also in \emph{log scale}, corresponds to $\ln{(\tilde{L}(c_0,c_1)-L_0)}$ with $c_0,c_1 \in [-0.1, 0.1]$, where $L_0$ is the smallest loss encountered overall.}
	\label{fig:mnist_complex-loss_manifold}
\end{figure}
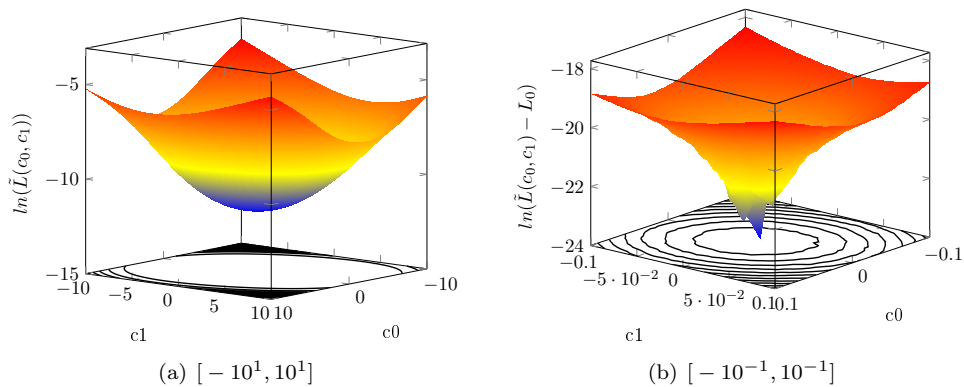

In Figure~\ref{fig:mnist_complex-loss_manifold} we give the resulting subspace manifold of the loss landscape for the sigmoid hidden activations function and the softmax cross-entropy loss. For the network using a 100 hidden nodes we have found a quasi-minimum during optimization that achieves a perfect accuracy of $1.0$ on the testset. Its accuracy on the training set is slightly lower with $0.975$.

In Fig.~\ref{fig:mnist_complex-loss_manifold-sigmoid-absolute} we look at the overall shape of the loss landscape in logarithmic representation at a length scale of $10$. We find a very smooth funnel. On a smaller scale of length $0.1$ in Fig.~\ref{fig:mnist_complex-loss_manifold-sigmoid-global}, where we look again at the logarithmic difference to the smallest loss encountered overall, we notice that there are multiple local minima close to the funnel's bottom. Not shown is the sampled domain of length scale $0.01$ where we then would see a heavily corrugated landscape.

This suggests that we see a similar landscape as with the linear model before where the bottom of the funnel is corrugated and that there is a plethora of states with small loss values.

\subsection{Application: \Glsentryname{glos:EQN} Method}\label{sec:application-ensemble_quasi_newton}
Returning to the analysis of the linear model, we have seen there is a large degree of anisotropy in the funnel observed in the loss manifold for the MNIST dataset. This results in significant variation in the projections of the gradient into different subspaces and is precisely the situation which motivated the development of the ensemble quasi-Newton scheme of Section 2.   

%This is a typical situation also encountered in \glostext{SD} methods in solving linear systems of equations where the matrix has a large condition number. A well-known solution is to use preconditioning.

%For sampling, the invariance of the sampled distribution under preconditioning needs to be ensured. One property that ascertains it is when the scheme is affine invariant: the distribution $\pi(\theta^{(n)}, p^{(n)})$ is invariant under affine transformations of the positions, $\pi(A \theta^{(n)} + v, p^{(n)})$. This depends on the choice of the preconditioner.

TATi facilitates the implementation of a scheme such as the ensemble quasi-Newton method (see Section 2) by it's multiple walker framework. In Appendix \ref{sec:eqn}, we describe the TATi implementation. Next, we investigate the method's qualities in a simple, well understood model before returning to the MNIST dataset.

\subsubsection{Gaussian Model}

A prototypical setting whose analytical properties are well-known is given by sampling from the Gaussian model,
\begin{equation}\label{math:gaussian_mixture_model}
	\exp{( -w^T C w)}
\end{equation}
with a covariance matrix $C \in \reel^{n \times n}$. Here, we naturally encounter directions that are \enquote{slow} to sample, identified by large eigenvalues in $C$. In fact, when we (only) look at the covariance structure of the MNIST loss manifold, we replace it by an effective Gaussian model of that particular covariance matrix. %Hence, it is only logical to study such a model as a first step.

Transfering this model to the setting of sampling loss manifolds of neural networks is straight-forward:
We use the mean squared loss $l_\theta \bigl (f (\theta, x_i), y_i \bigr ) = \bigl ( f(\theta, x_i) - y_i \bigr )^2$ with the network's prediction $f(\theta, x_i)$. Then, inserting into \eqref{math:lossfunction} in the case of $n$-dimensional input data $x_i \in \reel^n$, single-dimensional output $y_i \in \reel$, and a single-layer perceptron, i.e., $f_\theta(x) = w \cdot x + b$ with the parameters $\theta = \{w,b\}$ in the form of weights $w \in \reel^n$ and of a bias $b \in \reel$, we obtain

\begin{equation*}
	L(\theta, D) = \sum_i (w \cdot x_i + b - y_i)^2.
\end{equation*}

Setting the bias $b$ and all outputs $y_i$ to zero, we get

\begin{equation*}
	L(\theta, D) = \sum_i (w \cdot x_i)^2 = \sum_i \sum^n_{l,m=1} w_l (x_{l,i} x_{m,i}) w_m = \sum^n_{l,m=1} w_l \Bigl ( \sum_i x_{l,i} x_{m,i} \Bigr ) w_m.
\end{equation*}

Note that we sample from the canonical Gibbs distribution $\exp{\bigl (-\beta L(\theta, D) \bigr )}$. % with respect to \eqref{math:lossfunction}.
Therefore, the dataset needs to consist of rank-1 factors $x_i$ that represent the chosen covariance matrix with components $C_{lm} = \sum_i x_{l,i} x_{m,i}$ in order to match this with \eqref{math:gaussian_mixture_model}. These factors can be obtained for example through an eigendecomposition $C = V \Lambda V^T$ as the eigenvectors $V_i$ times the square root of their associated eigenvalue $\Lambda_{i,i}$.
Naturally, any other (even non-orthogonal) decomposition into rank-1 factors would be admissible, too.

In order to produce random covariance matrices $C$ of a certain structure, we
resort to the following approach: We generate a random symmetric matrix, compute its eigendecomposition and modify the diagonal matrix $D$ to consist of values picked from an equidistant spacing of the interval $[1,100]$, where endpoints are included. This way we obtain a set of orthogonal vectors pointing uniformly randomly in $\reel^n$, see \cite{Mezzadri2007}, and we make sure to generate both slow (eigenvalues close to 100) and fast (eigenvalues close to 1) directions. 

Having generated a random covariance matrix $C$ for dimensions $n \in \{2,4,8,16,32,64,128\}$ and having created the resulting dataset as its rank-1 factors, we sample the mean squared loss manifold of the single-layer perceptron using the BAOAB sampler. We use \num{50000} steps with a time step size $\varepsilon = 0.125$, inverse temperature $\beta = 1$, friction constant $\gamma = 1$, and covariance blending factor of $\eta = 10$.

We measure the exploration speed in the Gaussian model by looking at the \glostext{IAT} $\tau$ per random direction, that we know from the random matrix' eigendecomposition, by projecting it onto each eigenvector and measure the \glostext{IAT} using the package \emph{acor}~(\cite{acor}).

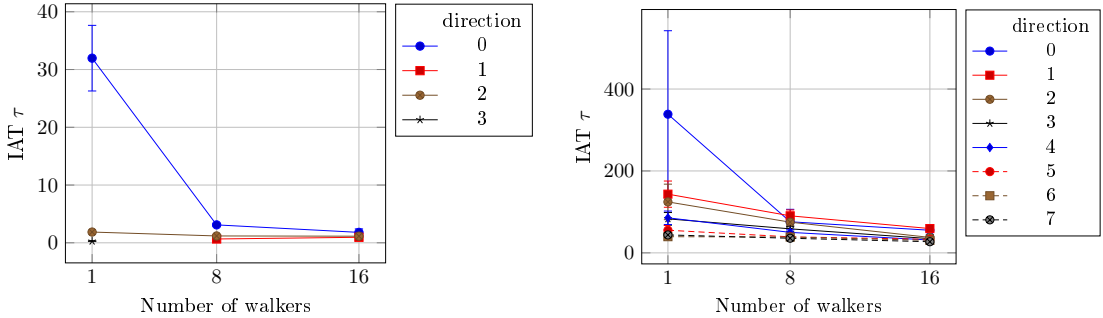
\begin{figure}[htbp]
	\resizebox{0.49\textwidth}{!}{
		\def\mydim{4}%
		\tikzsetnextfilename{autocorrelation_time-covariance_blending-collapse_walkers-dim_4}
		\ifdefined\myexperimentname\empty\else\def\myexperimentname{Gaussian Mixture}\fi%
\ifdefined\mypath\empty\else\def\mypath{\myexperimentspath/14_EnsembleQuasiNewton_GaussianMixture/data/Gaussian_Mixture}\fi%

%special column names
\ifdefined\myxcolumnname\empty\else\def\myxcolumnname{number_walkers}\fi%
\ifdefined\myycolumnname\empty\else\def\myycolumnname{tau}\fi%
%parameters
\ifdefined\mycovarianceblending\empty\else\def\mycovarianceblending{10.}\fi% %10
\ifdefined\mycollapsewalkers\empty\else\def\mycollapsewalkers{0}\fi% %0
\ifdefined\mydim\empty\else\def\mydim{32}\fi% %0

\begin{tikzpicture}[]
\begin{axis}[
	width=7cm,height=6cm,
	xlabel = {Number of walkers},
	ylabel = {IAT $\tau$},
	xtick={1,8,16},
	legend pos=outer north east,
	grid=major,
%	each nth point=1,
]
\addlegendimage{empty legend}
\addlegendentry{direction}
%\foreach \mydim in {2,4,8}{%
%\xdef\mycovarianceblending{\fpeval{10/\mydim}}%
\pgfmathparse{min(int(\mydim-1),7)}%
\xdef\mylimit{\pgfmathresult}%
%\node(debug) at (axis cs:8,10) {weights \mylimit};
\expandafter\pgfplotsinvokeforeach\expandafter{0,...,\mylimit}{%
 \addplot+[
	restrict expr to domain={\thisrow{weight}}{#1:#1},
	unbounded coords=discard,
	error bars/.cd,
        y dir=both,
        y explicit
] table[x=\myxcolumnname,y=avg_\myycolumnname, y error expr={sqrt(abs(\thisrow{var_\myycolumnname}-\thisrow{avg_\myycolumnname}^2))}, col sep=comma]{\mypath/autocorrelation_time-collapse_walkers_\mycollapsewalkers-covariance_blending_\mycovarianceblending-dimension_\mydim.csv};
\expandafter\addlegendentry\expandafter{#1}
}%
%}%
\end{axis}
\end{tikzpicture}
	}
%	\resizebox{0.32\textwidth}{!}{
%		\def\mydim{8}%
%		\tikzsetnextfilename{autocorrelation_time-covariance_blending-collapse_walkers-dim_8}
%		\input{./figures/14_EnsembleQuasiNewton_GaussianMixture/autocorrelation_time-covariance_blending-collapse_walkers}
%	}
	\resizebox{0.49\textwidth}{!}{
		\def\mydim{64}%
		\tikzsetnextfilename{autocorrelation_time-covariance_blending-collapse_walkers-dim_64}
		\ifdefined\myexperimentname\empty\else\def\myexperimentname{Gaussian Mixture}\fi%
\ifdefined\mypath\empty\else\def\mypath{\myexperimentspath/14_EnsembleQuasiNewton_GaussianMixture/data/Gaussian_Mixture}\fi%

%special column names
\ifdefined\myxcolumnname\empty\else\def\myxcolumnname{number_walkers}\fi%
\ifdefined\myycolumnname\empty\else\def\myycolumnname{tau}\fi%
%parameters
\ifdefined\mycovarianceblending\empty\else\def\mycovarianceblending{10.}\fi% %10
\ifdefined\mycollapsewalkers\empty\else\def\mycollapsewalkers{0}\fi% %0
\ifdefined\mydim\empty\else\def\mydim{32}\fi% %0

\begin{tikzpicture}[]
\begin{axis}[
	width=7cm,height=6cm,
	xlabel = {Number of walkers},
	ylabel = {IAT $\tau$},
	xtick={1,8,16},
	legend pos=outer north east,
	grid=major,
%	each nth point=1,
]
\addlegendimage{empty legend}
\addlegendentry{direction}
%\foreach \mydim in {2,4,8}{%
%\xdef\mycovarianceblending{\fpeval{10/\mydim}}%
\pgfmathparse{min(int(\mydim-1),7)}%
\xdef\mylimit{\pgfmathresult}%
%\node(debug) at (axis cs:8,10) {weights \mylimit};
\expandafter\pgfplotsinvokeforeach\expandafter{0,...,\mylimit}{%
 \addplot+[
	restrict expr to domain={\thisrow{weight}}{#1:#1},
	unbounded coords=discard,
	error bars/.cd,
        y dir=both,
        y explicit
] table[x=\myxcolumnname,y=avg_\myycolumnname, y error expr={sqrt(abs(\thisrow{var_\myycolumnname}-\thisrow{avg_\myycolumnname}^2))}, col sep=comma]{\mypath/autocorrelation_time-collapse_walkers_\mycollapsewalkers-covariance_blending_\mycovarianceblending-dimension_\mydim.csv};
\expandafter\addlegendentry\expandafter{#1}
}%
%}%
\end{axis}
\end{tikzpicture}
	}
	\caption{\protect\glosfirst{IAT} $\varepsilon$ over the number of walkers for Gaussian models of dimensions $n \in \{4,64\}$ where only up to the first 8 directions are shown. Five Sampling runs have been completed for each random covariance matrix and average \protect\glostext{IAT} and standard deviation is shown. With a single walker standard BAOAB sampling is employed, with multiple walkers the \protect\glostext{EQN} method is used.}
	\label{fig:gaussian_mixture-iat}
\end{figure}

In Figure~\ref{fig:gaussian_mixture-iat} we see that using the \glostext{EQN} scheme with 8 or 16 walkers significantly improves the Integrated Autocorrelation Time $\varepsilon$ for the slow directions. We note that the fast directions are unaffected. We remind the reader that each walker samples its own trajectory. Hence, we generate up to 16 trajectories in parallel and do so ten times more efficiently because of the reduction in the IATs.

\subsubsection{MNIST}

We now turn to the MNIST dataset again for a real-world application of the \glostext{EQN} method.
We constrain the training dataset to two classes, namely the digits 7 and 9. This results in \num{11169} training dataset items for this two-class problem. We employ the same single-layer perceptron as before.

For the \glostext{IAT} computation we have extracted distinct covariance matrices $C_i = \cov(\theta_i^{(n)},\theta_i^{(n)})$ per walker, obtained the eigendecomposition $C = V \Lambda V^T$, extracted the eigenvectors $V_j$ of the 20 dominant eigenvalues, and projected the walker's trajectory onto these, $\pi_{i,j} (n) = \Theta_i^{(n)} \cdot V_j$. Finally, we calculated the \glostext{IAT} of each $\pi_{i,j} (n)$ using the \module{acor} package and averaged over all walkers $i$.

Using walker-individual covariance matrices does not generally change results compared to a single covariance matrix obtained from averaging the trajectory over all walkers; however, it makes them more stable. Because of the high dimensionality of  the parameter space $\reel^{\num{\fpeval{784*2}}}$ already small perturbations may cause vectors to become orthogonal\footnote{The more dimensions a space has, the more likely it becomes for two random vectors to be orthogonal to each other, see also \cite[\citesection~7.1]{Li2017}.}. This phenomena is explained by the \enquote{Concentration of Measure}, see \cite{LedouxMichel2001Tcom}. The eigenvalue spectra themselves are stable for each walker, bounded in deviation by the Bauer-Fike theorem.
Note further that the preconditioning matrix $B_i^{(n)}$ is also uniquely defined for each walker.

Note that because of computations necessary for the additional precondition matrix, there is a slight walker-dependent overhead compare to the case of single walker: For 8 walkers we measured about 40\% overhead per walker, for 16 walkers we obtained 50\%.

% show IAT improvements on the MNIST dataset
%

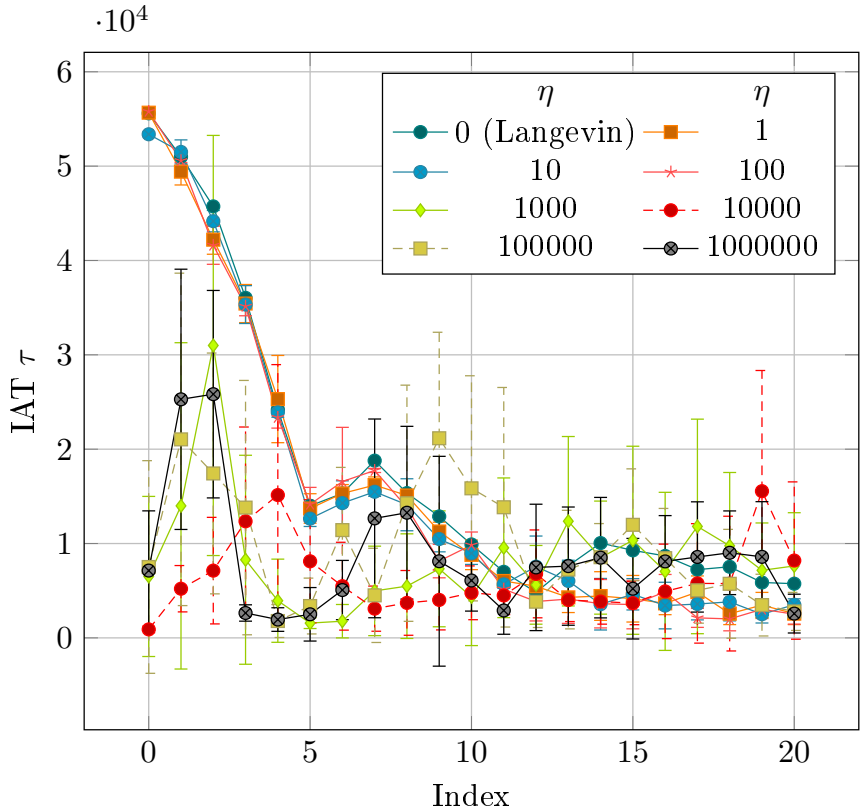
\begin{figure}[htbp]
	\centering
	\resizebox{0.8\textwidth}{!}{
		\def\myinversetemperature{10} %
		\def\mynumberwalkers{8}%
		\def\mybatchsize{550}%
		\tikzsetnextfilename{MNIST-autocorrelation-individual_time-covariance_blending-linear_scale}
		\ifdefined\myexperimentname\empty\else\def\myexperimentname{Parameter study eta - Fixing biases}\fi%
\ifdefined\mypath\empty\else\def\mypath{\myexperimentspath/19_EQN-MNIST-two_classes/data/Parameter_study_eta_-_Fixing_biases}\fi%

%special column names
\ifdefined\myxcolumnname\empty\else\def\myxcolumnname{weight}\fi%
\ifdefined\myycolumnname\empty\else\def\myycolumnname{mean}\fi%
%parameters
\ifdefined\myeverynth\empty\else\def\myeverynth{100}\fi% %100
\ifdefined\myinversetemperature\empty\else\def\myinversetemperature{10}\fi% %10
%\ifdefined\mycovarianceblending\empty\else\def\mycovarianceblending{10}\fi% %10
\ifdefined\mybatchsize\empty\else\def\mybatchsize{550}\fi% %550
\ifdefined\mynumberwalkers\empty\else\def\mynumberwalkers{8}\fi% %8

\begin{tikzpicture}[]
\begin{axis}[
	width=10cm,height=9cm,
	xlabel = {Index},
	ylabel = {IAT $\tau$},
	legend pos=north east,
	legend columns=2,
	cycle list name=exotic,
	grid=major,
%	each nth point=1,
	restrict x to domain={0:20},
]
\addlegendimage{empty legend}
\addlegendentry{$\eta$}
\addlegendimage{empty legend}
\addlegendentry{$\eta$}
% 
% \addplot+[no marks, thick] table[x=\myxcolumnname,y=avg_tau, col sep=comma]{\myexperimentspath/19_EQN-MNIST-two_classes/data/Sampling_-_EQN_TATi_Single/autocorrelation-individual_time-batch_size_\mybatchsize-covariance_blending_0.csv};
% \addlegendentryexpanded{1}

\foreach \mycovarianceblending/\mycovarianceblendingname in {0/0 (Langevin),1/1,10/10,100/100,1000/1000,10000/10000,100000/100000,1000000/1000000}{%
\addplot+[error bars/.cd, y dir=both, y explicit]
 table[x=\myxcolumnname,y expr={\thisrow{avg_tau}*\myeverynth}, y error expr={sqrt(\thisrow{var_tau}-\thisrow{avg_tau}^2)*\myeverynth}, col sep=comma]{\mypath/autocorrelation-individual_time-batch_size_\mybatchsize-covariance_blending_\mycovarianceblending-inverse_temperature_\myinversetemperature-number_walkers_\mynumberwalkers.csv};
\expandafter\addlegendentry\expandafter{\mycovarianceblendingname}
}%
\end{axis}
\end{tikzpicture}
	}
	\caption{\protect\glosfirst{IAT} $\tau$ for the first 20 covariance eigenvectors over the number of walkers for the MNIST two class problem. We observe a strong improvement of up to a factor of 4 for the slowest \protect\glosfirst{IAT} which relates to a similar increase in exploration speed, e.\,g., when using $\eta$ of \num{10000} (EQN)compared to $\eta = 0$ (standard Langevin).}
	\label{fig:mnist-iat}
\end{figure}

In Figure~\ref{fig:mnist-iat} we then look at the \glostext{IAT} over the first 20 covariance eigenvectors for various values of the covariance blending constant $\eta$. We used a fixed number of 8 walkers, a batch size of \num{550}, an inverse temperature constant of $\beta = 10$, the friction constant set to $\gamma = 10$ and a step size $\varepsilon = 0.125$ with the BAOAB sampler. All runs are started from an equilibrated position using \glostext{SGD} with batch size of \num{550} and a learning rate of $0.1$ for \num{5000} steps.  We recompute the covariance matrix after \num{10000} steps.
The trajectories were stored with only every 100th sampling step. Hence, the $\varepsilon$ values in the Figure have been rescaled appropriately.

As there is a scaling invariance with respect to the biases of the output layer due to the \python{argmax} function, we have fixed the biases for the sampling to the values obtained from a prior optimization. At the moment, the \glostext{EQN} implementation cannot deal with such invariances. They represent \enquote{flat valleys} in the loss landscape and the walkers will be pushed by the preconditioning along the valley in vain search for its bounds. Such a valley can be hypothesized from the spectrum of the covariance matrix, see Figure~\ref{fig:mnist_sampling-covariance_eigenvalues}, where the first eigenvalue with $2.1 \cdot 10^6$ is unusually high, and when the \glostext{EQN} does not effectively reduce the \glostext{IAT}s although the eigenvalue spectrum indicates it, c.\,f.~\cite[\citepage~281]{Matthews2018}.

We generally observe that especially the first five IAT values are dramatically reduced and see an improvement of up to a factor of 4.  The covariance blending can be chosen robustly, up to very high values. This indicates that the covariance matrix itself, despite $L < N$, is already positive definite.

% - compare (unpreconditioned) IATs with eigenvalue structure from covariance matrix. Mind a possible every_nth!
%
% - show preconditioned IATs for multiple values of eta.

\section{Conclusion}\label{sec:conclusions}
We have discussed comparisons of sampling strategies based primarily on stochastic differential equations, with a focus on issues such as sampling efficiency and accuracy.  All the methods discussed are implemented in the TATi software.   Relying on \TF{}, the TATi implementation efficiently runs in parallel and also on GPU-assisted hardware.

We have looked in our evaluation at the MNIST loss manifold for a single-layer perceptron and a multi-layer perceptron using softmax cross-entropy. We find that it resembles an anisotropic funnel on the large scale combined with many local minima at its bottom matching well the band of minima  bounded from below and exponentially decaying in density with higher loss values predicted in \cite{Choromanska2015}.
This motivated an ensemble method employing a number of so-called walkers to obtain a local approximation of the covariance that, when turned into a preconditioner, results in a significant improvement of the sampling speed.

The focus on posterior sampling using SDEs provides us with a starting point for a wide range of improvements.  For example we are currently exploring schemes based on simulated tempering~\cite{ST,ISST} and diffusion maps~(\cite{Coifman2006a,Chiavazzoa2017}); the latter allows simulation data obtained using e.g. Langevin dynamics to be distilled into a few collective variables which succinctly describe the progress of transition between neighboring local minima.   We are also exploring the use of the sampling paradigm as a tool to design neural networks with improved sparsity and generalizability.   

\appendix

\section*{Addendum}

\section{Virial Theorem}\label{sec:virial_theorem}
% show that virial theorem only holds for for potentials unbounded from above that raise at least polynomially
The virial is defined as $G = \sum^N_i \theta_i p_i$ with positions $\theta$ and momenta $p$. The virial theorem states that $\frac{d G}{dt} = 0$ which would imply through 
\begin{equation}\label{math:virial_theorem}
\frac{d G}{dt} = \sum^N_i p_i p_i + \sum^N_i \theta_i \frac{\partial L(\theta)}{\partial \theta_i}
\end{equation}
that the first term, twice the kinetic energy, equals the negative of the second.

Let us inspect the second term and look at its average over the whole domain using the Gibbs measure with a single degree of freedom  ($N=1$),
\begin{equation*}
	\frac{ \int_\R \theta \cdot \nabla L(\theta) \exp{(-\beta L(\theta))} d\theta }{ \int \exp{(-\beta L(\theta))} d\theta  },
\end{equation*}
where $\beta$ is the inverse temperature.

Let us ignore the denominator for the moment and integrate the nominator by parts. We obtain with the derivative $\tfrac{\partial}{\partial \theta} \exp{(-\beta L(\theta))} = - \beta \exp{(-\beta L(\theta))} \tfrac{\partial L(\theta)}{\partial \theta}$,
\begin{equation*}
	\int_\R \theta \cdot \nabla L(\theta) \exp{(-\beta L(\theta))} d\theta = \Bigl [ \theta \cdot \tfrac{-1}{\beta} \exp{(-\beta L(\theta))} \Bigr ] - \int_\R 1 \cdot \tfrac{-1}{\beta} \exp{(-\beta L(\theta))} d\theta.
\end{equation*}

If the boundary term vanishes, we obtain
\begin{equation*}
	\frac{ \int_\R \frac 1 \beta \exp{(-\beta L(\theta))} d\theta }{ \int_\R \exp{(-\beta L(\theta))} d\theta  } = \tfrac 1 \beta.
\end{equation*}
In other words, the average virial, the second term, would be identical to two times the average kinetic energy $\tfrac 1 {\beta}$, the first term in \eqref{math:virial_theorem}.

Hence, all that remains is to show that $\Bigl [ \theta \cdot \exp{(-\beta L(\theta))} \Bigr]^{\infty}_{-\infty} = 0$. Naturally, this holds if $\lim_{|\theta| \rightarrow \infty} L(\theta) \rightarrow \infty$ to the effect that $\exp{(-\beta L(\theta))} \rightarrow 0$ faster than $|\theta| \rightarrow \infty$, noting that the exponential increases faster than any polynomial.

In other words, the potential $L(\theta)$ needs to be unbounded and to increase faster than $|\theta|$ for the virial theorem to hold.

\subsection{Virial and MNIST}\label{sec:virial_theorem-softmax_cross_entropy}

If for the MNIST dataset, a single-layer perceptron with a softmax cross-entropy function is employed, then the virial theorem does not hold. 

The output of the single-layer perceptron is $f_i(\theta) = \sum_j W_{i,j} x_j + b_i$ with weight matrix $\theta$ and bias vector $b$, i.\,e. $\theta = (W,b)$.

%When for any component $i,j$ of the parameter matrix $\theta$, we have $\theta_{i,j} \rightarrow \infty$ of a single-layer perceptron , then we have for its output $f_i(\theta) \rightarrow \infty$ assuming at least one non-zero input $x_j$ in ${\cal D}$. The softmax function $p_i(f(\theta)) = \tfrac{ \exp{f_i(\theta)} }{ \sum_i \exp{f_i(\theta)} }$ bounds it from above to 1, i.\,e.~$p_i(f(\theta)) \rightarrow 1$ for $\theta_{i,j} \rightarrow \infty$. Therefore, $L(\theta) = - \sum_j y_j \log{p_j(\theta)}$ is bounded from above for the cross entropy.

As the cross entropy is $- \sum_i y_i \log{p_i(f(\theta))}$ and the softmax function is $p_i(f(\theta)) = \tfrac{ \exp{f_i(\theta)} }{ \sum_i \exp{f_i(\theta)} }$, we have $\exp{(\beta \sum_i y_i\log{p_i(f(\theta))})} = \prod_i p_i(f(\theta))^{\beta y_i}$.

Let us set all parameters components to zero except for a single weight component $W_{i,j}$ where at least for one item in the dataset we have $x_j \neq 0$.
Then we obtain $W_{i,j} (W_{i,j} x_j)^{\beta y_i}$ as the integrand for this data item $(x_i, y_i)$ and the boundary integral will not converge (to zero) in this case.

Note that this issue could be addressed by the addition of an $L_2$ regularization strategy.

% Therefore, $\beta y_i$ needs to be smaller than $-1$ for at least one $y_i$ for the potential to decrease sufficiently quickly.

%\bibliographystyle{alpha}
%\bibliography{sample}

% Curse of Dimensionality: fast and slow directions

\section{Ensemble Quasi-Newton implementation in TATi}\label{sec:eqn}
In Listing~\ref{lst:eqn_sampler} we provide a rapid prototype of the algorithm using \glostext{TATi}'s \module{simulation} module.

\begin{lstlisting}[
	language=Python,
	basicstyle=\tiny,
	caption={Example implementation of \protect\glostext{EQN} using \protect\glostext{TATi}'s \module{simulation} module. We require for the number of walkers $L > 1$. We have skipped (\ldots) the details of instantiation of the interface class \python{tati} setting the options.},
	label={lst:eqn_sampler},
]
import math
import numpy as np
import tensorflow
import TATi.simulation as tati

nn = tati(
    ....
)
options = nn.get_options()

def baoab_update_step(nn, momenta, old_gradients, preconditioner, step_width, beta, gamma, walker_index):
    def B(step_width, gradients):
        nonlocal momenta
        momenta -= .5 * step_width * np.dot( np.transpose(preconditioner), gradients)

    def A(step_width, momenta):
        nn.parameters[walker_index] += .5 * step_width * preconditioner.dot(momenta)

    def O(step_width, beta, gamma):
        nonlocal momenta
        alpha = math.exp(-gamma * step_width)
        momenta = alpha * momenta + \
                  math.sqrt((1. - math.pow(alpha, 2.)) / beta) * np.random.standard_normal(momenta.shape)

    B(step_width, old_gradients)
    A(step_width, momenta)
    O(step_width, beta, gamma)
    A(step_width, momenta)
    gradients = nn.gradients(walker_index=walker_index)
    B(step_width, gradients)

    return gradients, momenta

preconditioner = [np.identity((nn.num_parameters())) for i in range(options.number_walkers)]
normalization = 1./(float(nn.num_walkers()) - 1.)

def update_preconditioner():
    for walker_index in range(nn.num_walkers()):
        means = normalization*np.add.reduce([nn.parameters[i] for i in range(nn.num_walkers()) if i != walker_index])
        covariance = normalization*np.add.reduce(\
            [np.outer(nn.parameters[i] - means, nn.parameters[i] - means) \
            for i in range(nn.num_walkers()) if i != walker_index])
        preconditioner[walker_index] = np.linalg.cholesky( \
            options.covariance_blending * covariance + np.identity((nn.num_parameters())))

momenta = [np.zeros((nn.num_parameters())) for i in range(options.number_walkers)]
old_gradients = [np.zeros((nn.num_parameters())) for i in range(options.number_walkers)]

def perform_step():
    for walker_index in range(nn.num_walkers()):
        old_gradients[walker_index], momenta[walker_index] = baoab_update_step(
            nn, momenta[walker_index], old_gradients[walker_index],
            preconditioner=preconditioner[walker_index],
            step_width=options.step_width,
            beta=options.inverse_temperature, gamma=options.friction_constant, walker_index=walker_index)

for step in range(options.max_steps):
    if (step) % options.covariance_after_steps:
        update_preconditioner()
    perform_step()
\end{lstlisting}

The full implementation with \TF{} in \glostext{TATi} requires special care with the conditionals for the infrequent updates of $B_i$, needs to copy the network parameters to avoid changes within the parallel execution, and needs to compute the covariance matrices.
All implemented samplers have been adapted in a similar way as in \eqref{math:baoab-preconditioned} to allow for preconditioning.
For performance reasons the computation of $B_i^{(n)}$ could be done entirely through rank-1 updates. 
%However, at the time of writing this is not yet implemented.

\section*{Acknowledgements}\label{sec:acknowledgements}

The authors acknowledge funding through a Rutherford fellowship from the Alan Turing Institute in London (R-SIS-003, R-RUT-001) and EPSRC grant no. EP/P006175/1 (Data Driven Coarse Graining using Space-Time Diffusion Maps, B. Leimkuhler PI), as well as B. Leimkuhler's Turing Fellowship (The Alan Turing Institute is supported by EPSRC grant EP/N510129/1).   The computing resources were provided by a Microsoft Azure Sponsorship award (MS-AZR-0143P).

In relation to the development of the TATi software, the authors would also like to express their gratitude to the Research Software Engineering Team of the Alan Turing Institute led by James Hetherington and especially to Martin O'Reilly for support in using the Azure platform.
Charles Matthews provided valuable comments on a preliminary draft of the article, leading to a much improved exposition.
%\ifloadtikztodos
%  \listoftodos
%\else
%  \todos
%\fi%

\bibliographystyle{plainnat}	% either alpha for [Sta99] or plain for [01]  or unsrt
\bibliography{small_library}
%\printbibliography

\printglossaries

\end{document}
We have implemented Gradient Descent with a Barzilei-Borwein step width choice (\glsentrytext{glos:BBGD}), see \cite{Tan2016}, using the full gradient information. We use 3000 \glsentrytext{glos:BBGD} steps with an initial learning rate of 0.05. This advanced optimization method is not required for the sampling's starting point, there \glostext{SGD} would suffice. However, it yields an improved or